%% file: involutions7.tex
\documentclass[a4paper,12pt,twoside]{article}
\usepackage{amssymb}
\usepackage{mathtools}
\usepackage[final]{lpretex}
\usepackage{title}
\usepackage{a4}
\usepackage{amscd}
\input xy
\xyoption{all}

\EndOfDump

\def\DHpartformat#1{(#1) }
\def\DHrefpart#1{(\DHRefpart{#1})}
\LBsetstyleof{partno}{arabic}

\DocVersion[Changed to LaTeX]{1}

\let\define\def

\def\A {{\mathbb A}}  \def\C {{\mathbb C}}
  \def\F {{\mathbb F}}

\def\GG {{\mathbb G}}   
  
\def\N {{\mathbb N}}  \def\P {{\mathbb P}} 
\def\Q {{\mathbb Q}} \def\R {{\mathbb R}}

\def\Z {{\mathbb Z}} 

\define \n {\mathbb N}
\define \z {\mathbb Z}
\define \q {\mathbb Q}
\define \PP {\mathbb P}

\def\sA {{\Cal A}}  
 \def\sE {{\Cal E}} \def\sF {{\Cal F}}
\def\sG {{\Cal G}} 
\def\sH {{\frak H}}

 \def\sN {{\Cal N}} \def\sO {{\Cal O}}
 
 \def\sT {{\Cal T}} \def\sU {{\Cal U}}

\define \cN {\Cal N}
\define \cf {\Cal F}
\define \cg {\Cal G}
\define \cE {\Cal E}
\define \ce {\Cal E}
\define \cc {\Cal C}
\define \cV {\Cal V}
\define \cA {\Cal A}
\define \cK {\Cal K}
\define \cO {\Cal O}
\define \cF {\Cal F}
\define \cn {\Cal N}
\define \cI {\Cal I}
\define \sP {\Cal P}
\define \sZ {\Cal Z}
\def\tA {\widetilde{\Cal A}}

\def\a {\alpha} \def\b {\beta}  \def \d {\delta} 
   
\def\s {\sigma} \def\t {\theta}

\define \x {\xi}
\define \y {\eta}
\define \G {\Gamma}
\define \r {\rho}
\define \w {\omega}

\def \tY {\widetilde Y}

\def \trho {\widetilde {\rho}}

\def \tp {\widetilde{\mathbb P}}
\define \tH {\widetilde H}
\define \tG {\widetilde{\Gamma}}
\define \tW {\widetilde W}
\define \tF {\widetilde F}
\define \tm {\widetilde m}
\define \St {\widetilde S}
\define \Xt {\widetilde X}
\define \tS {\widetilde S}
\define \tpsi {\widetilde \psi}
\define \tL {\widetilde L}
\define \tE {\widetilde E}
\define \tl {\widetilde l}
\define \tA {\widetilde A}
\define \tom {\widetilde\omega}
\define \tT {\widetilde T}
\define \tB {\widetilde B}
\define \tf {\widetilde f}
\define \tsA {\widetilde{\sA}}

\define \tM {\widetilde M}
\define \tphi {\widetilde{\phi}}
\define \trho {\widetilde{\rho}}
\define \tR {\widetilde R}
\define \tp {\widetilde p}
\define \tq {\widetilde q}
\define \tc {\widetilde c}
\define \tsF {\widetilde {\sF}}
\define \tsN {\widetilde {\sN}}
\define \tsU {\widetilde {\sU}}
\define \th {\widetilde h}
\define \tsZ {\widetilde {\sZ}}

\def\pd {\partial}

\def \Dx1 {\frac{\pd}{{\pd} x_1}}
\def \Dy1 {\frac{\pd}{{\pd} y_1}}
\def \Dz1 {\frac{\pd}{{\pd} z_1}}
\def \Dx2 {\frac{\pd}{{\pd} x_2}}
\def \Dy2 {\frac{\pd}{{\pd} y_2}}
\def \Dz2 {\frac{\pd}{{\pd} z_2}}

\def\q {\quad} 

\def\mapdiagr#1{\Big\searrow\rlap{$\raise 5pt\vbox{{\hbox{$\mkern -15mu\scriptstyle#1$}}}$}}   

\def\mapdiagl#1{\llap{$\raise 5pt\vbox{{\hbox{$\scriptstyle#1\mkern
-15mu$}}}$}\Big\swarrow}              

\def\Mapdiagr#1{\nearrow\rlap{$\lower 5pt\vbox{{\hbox{$\mkern
-15mu\scriptstyle#1$}}}$}} 

\def\Mapdiagl#1{\llap{$\lower 5pt\vbox{{\hbox{$\scriptstyle#1\mkern
-15mu$}}}$}\searrow} 

\def\Mapswr#1{\swarrow\rlap{$\lower 5pt\vbox{{\hbox{$\mkern
-15mu\scriptstyle#1$}}}$}}              

\def\Mapnwl#1{\nwarrow\rlap{$\lower 5pt\vbox{{\hbox{$\mkern
-15mu\scriptstyle#1$}}}$}}

\define \Rhook {\hookrightarrow}

\def \half {\raise1pt\hbox{$\scriptstyle
        \frac{1}{2}\displaystyle$}}

\def \x{{\sl X}\llap{$\mkern -2mu {\scriptstyle -}$}}
\def \Bl {\operatorname{Bl}}

\define \Kod {\operatorname{Kod}}
\define \dimension {\operatorname{dim}}
\define \codim {\operatorname{codim}}
\define \contr {\operatorname{contr}}
\define \rk {\operatorname{rank}}
\define \im {\operatorname{im}}
\define \Mor {\operatorname{Mor}}
\define \Cl {\operatorname{Cl}}
\define \Hilb {\operatorname{Hilb}}
\define \degree {\operatorname{deg}}
\define \mult {\operatorname{mult}}
\define \Aut {\operatorname{Aut}}
\define \NS {\operatorname{NS}}
\define \Gal {\operatorname{Gal}}
\define \ch {\operatorname{char}}
\define \Jac {\operatorname{Jac}}
\define \Km {\operatorname{Km}}
\define \Sec {\operatorname{Sec}}
\define \Stab {\operatorname{Stab}}
\define \Br {\operatorname{Br}}
\define \inv {\operatorname{inv}}
\define \tr {\operatorname{tr}}
\define \Frob {\operatorname{Frob}}
\define \Symn {\operatorname{Sym}^n}
\define \Ev {\sE^\vee}
\define \ordp {\operatorname{ord}_p}
\define \Supp {\operatorname{Supp}}
\define \Ann {\operatorname{Ann}}
\define \disc {\operatorname{disc}}
\define \Lie {\operatorname{Lie}}
\define \embdim {\operatorname{embdim}}

\def \Fix{\operatorname{Fix}}

\def\Tr{\operatorname{Tr}}

\def\diam{\operatorname{diam}}
\def\Tub{\operatorname{Tub}}

\def\Bir{\operatorname{Bir}}
\def\Ax{\operatorname{Ax}}

\def\eps{\epsilon}

\def\bark{\overline{k}}

\def\hod#1#2#3#4{\ensuremath{\if#30 H^{#2}({#1},{\cal O}_{#1}) \else 
 H^{#2}(#1,\Omega^{#3}\if\relax{#4}\relax_{#1}\else _{#1/#4}\fi)\fi}}

\begin{document}
\title[Effectivity questions in the Cremona group]
{Some effectivity questions for plane Cremona transformations
in the context of symmetric-key cryptography}
\author{N. I. Shepherd-Barron}
\address{Mathematics Department\\
King's College London\\
Strand\\
London WC2R 2LS\\
U.K.}
\email{Nicholas.Shepherd-Barron@kcl.ac.uk}

\maketitle
\begin{section}{Introduction}

The $2$--dimensional Cremona group $Cr_2(k)$ is the group of $k$-birational automorphisms
of the projective plane $\P^2_k$ over a field $k$. As such, it is an object
of algebraic geometry, but it is also of interest from the
viewpoints of dynamics [M1] and group theory, including geometric
group theory [CL]. 
This latter paper introduces and uses a certain
infinite-dimensional hyperbolic space $\sH=\sH_k=\sH(\P^2_k)$
on which $Cr_2(k)$ acts as a group of isometries
and makes it clear that this action is an important tool for studying
both $Cr_2(k)$ and its individual elements.
However, given an algebraic description of a Cremona
transformation $g$ in terms of explicit rational functions, it is
not always clear how to calculate anything about $g$
that is relevant to any of these frameworks.
For example, there is no known effective procedure for determining
the \emph{translation length} $L(g)=L(g_*)$ of the isometry 
$g_*$ of $\sH$ that is associated to $g$.

One of our goals here is to make matters more nearly effective in
the very particular case where $g$ is 
a \emph{special quadratic transformation}.
That is, $g=\s\a$, where $\s$ is the standard quadratic 
transformation given in terms of homogeneous co-ordinates
$x,y,z$ by $\s:(x,y,z){\mapstochar\dashrightarrow} (yz,xz,xy)$ 
and $\a$ is a linear involution.
For this special class of Cremona transformations
we give (Theorem \ref {Trrr} below)
a simple and explicit condition, in terms
of a constant number of field operations
(addition and multiplication in the ground field)
to ensure that $g$ is hyperbolic: it is enough to
take the three vertices $P,Q,R$
of the triangle $\Delta$ given by $xyz=0$,
compute the $12$ points $g^i(P)$ etc. for $1\le i\le 4$
and verify that they are distinct and that none of them lie on $\Delta$. 
In comparison, Bedford and Kim [BK, Theorem 3.2] 
give an exact formula
for $L(g)$ when $\a$ is any linear transformation
that relies on an unbounded quantity of information
about the map $g$.
The upper bound $L(h)\le\log D$ 
for a Cremona transformation $h$ of degree $D$
is well known [CL]; these bounds together
give an arbitrarily fine estimate for $L(g)$.
These bounds will also determine the type of $g$; that is,
whether $g$ is elliptic, parabolic or hyperbolic (the word
loxodromic is also used for this last class).
This is motivated by the fact, explained in Section \ref{crypto},
that certain contemporary encryption algorithms are
Cremona transformations and that estimates 
such as those proved here can be seen as a speedy check
that the key being used is not obviously weak.

Say that points $x_1,x_2,...,x_n$ in $\P^2$
are \emph{in general position} if they are all
distinct and none of them, except for those that happen to
equal one of $P,Q,R$, lies on $\Delta$.
For any natural number $n$, let $w_n$ denote the
reduced word in $\s,\a$ that has length $n$ and begins
with $\a$ when reading from right to left. So, for example,
$w_0=1,w_1=\a,w_2=\s\a$. Set $P_n=w_n(P)$, $Q_n=w_n(Q)$
and $R_n=w_n(R)$. So $P_0=P$, etc.
We say that $\a$, or $g$, is 
\emph{in $(p,q,r)$-general position}
if $P_0,...,P_{p-1},Q_0,...,Q_{q-1},R_0,...,R_{r-1}$
are in general position. We abbreviate $(r,r,r)$-general
position to $r$-general position.
Of course, $r$-general position implies $s$-general position
for any $s\le r$.

\begin{theorem}\label{theorem 1} 
\part \label{T444} (= Theorem \ref{4.21}) Suppose that $g$
is in $(p,q,r)$-general position
and that $1/p+1/q+1/r<1$.
Then $g$ is hyperbolic.
If also $p\le q\le r$ then
$L(g)\ge\log(2-3.2^{-p/2})$.

\part \label{Trrr} (= Theorem \ref{theorem a})
Assume that $0<\eps<1/3$.
Then, if $p\ge 12/\eps$ and
the $3p$ points $P_0,...,P_{p-1},...,R_{p-1}$
are in general position, 
$g$ is hyperbolic and
$$\log 2-\epsilon< L(g)\le \log 2.$$
\noproof
\end{theorem}

We can also give some analogous sufficient conditions
that permit the precise determination of the type
of a special quadratic transformation
$g$ and its length. Here is an example.

Assume that $g$ is in $(p,q,r)$-general position
and that
$P_p=P_{p-1}$, $Q_q=Q_{q-1}$ and $R_r=R_{r-1}$.
Let $\lambda_{Lehmer}$ denote \emph{Lehmer's number},
the smallest known algebraic integer $\lambda$
such that $\vert\lambda\vert>1$ and every other conjugate
$\lambda'$ of $\lambda$ has $\vert\lambda'\vert\le 1$
and $\vert\lambda'\vert=1$ for at least one $\lambda'$.
(These are the \emph{Salem numbers}.)

\begin{theorem}\label{theorem 2} (= Theorem \ref{Coxeter})
\part
The transformation $g$ is biregular on the blow-up
of $\P^2$ at these $p+q+r$ points and is a Coxeter
element in a Weyl group of type $T_{p,q,r}$.

\part If $1/p+1/q+1/r>1$ then $g$ is elliptic
and its order is the corresponding Coxeter number.
In particular, if $p=2$, $q=3$ and $r=5$ then $g$ has order $30$.

\part If $1/p+1/q+1/r=1$ then $g$ is parabolic.

\part If $1/p+1/q+1/r<1$ then $g$ is hyperbolic and $L(g)$
is the logarithm of a Salem number of norm $1$.
In particular, if $p=2$, $q=3$ and $r=7$
then $L(g)=\log\lambda_{Lehmer}$.
\noproof
\end{theorem}

These results, and their proofs, can be summarized by saying that
there is a Coxeter--Dynkin diagram
associated to the problem and
if $P_0,...,R_{r-1}$ are in general position then
the diagram contains the standard tree
$T_{p,q,r}$. Since $T_{2,3,5}=E_8$
the number $30$ appears as
the Coxeter number of $E_8$.
This is a particular instantiation of the very old idea
of relating groups of Cremona transformations
to Coxeter--Dynkin diagrams and
Weyl groups, although usually the diagrams that have arisen
are of type $T_{2,3,r}$, and also of Steinberg's idea
[S] of describing Coxeter elements as a product of two involutions
via a description
of the Coxeter diagram as a bipartite graph.
This viewpoint has also been exploited by
McMullen [M1], whose calculations have inspired some of 
those that appear here, and
Blanc and Cantat [BC], who use
an infinite group $W_\infty$ that is
something like a Coxeter group of type $E_\infty$
to prove that, for any hyperbolic Cremona transformation
$h$, the \emph{spectral radius} or \emph{dynamical degree}
$\lambda(h)=\exp L(h)$
lies in the closure of the set $\sT$ of Salem numbers and
$\lambda(h)\ge\lambda_{Lehmer}$.

The other main point of the paper is to extend the results
of [CL] that concern the \emph{rigidity} and \emph{tightness} of
elements, or conjugacy classes, in the Cremona group. 
(The definitions of these properties
are recalled later, at the start of Section \ref{7}.)
Some of these results also complement a recent paper [L] by Lonjou,
who exhibits, over any field $k$, an explicit Cremona transformation
some power of which generates a proper normal subgroup of $Cr_2(k)$.
She also points out an error in an earlier version of this paper;
the mistake lay in overlooking the possibility that a rigid element
of $Cr_2(k)$ might normalize a $2$-dimensional additive subgroup
of $Cr_2(k)$. In consequence, the results of Section \ref{7}
below in characteristic $p$ that refer to normal subgroups
of $Cr_2(k)$
require the assumption that $k$ be algebraic
(that is, algebraic over its prime subfield).

Fix a plane Cremona transformation $g$ over a field $k$.

\begin{theorem} Assume that $g$ is hyperbolic.

\part (= Theorem \ref{lots}) $g$ is rigid.

\part \label{unit}(= Theorem \ref{4.13}) 
Suppose that $L(g)$ is not the logarithm of a 
quadratic unit; if $\ch k=p>0$ assume also
that $k$ is algebraic and
that $L(g)$ is not an integral multiple of $\log p$. 
Then some power of $g$ is tight.
If also $n$ is sufficiently divisible,
then the normal closure $\langle\langle g^n\rangle\rangle$
does not contain $g$, so is a non-trivial normal subgroup of $Cr_2(k)$.
\noproof
\end{theorem}

In [CL] it was shown that if $g$ is a 
very general Cremona transformation
of any degree $\ge 2$, then some power of
$g$ is tight and that $g\notin\langle\langle g^n\rangle\rangle$
for sufficiently divisible $n$.
(If $g$ is very general then $L(g)=\log \deg(g)$.)

In particular, this applies to special quadratic transformations
$g$ if $k$ is algebraic and
$\ch k\ne 2$ (if $\ch k=2$ we must assume also that $L(g)\ne \log 2$), 
since it is an easy consequence of the other results that
we prove about them that they satisfy the hypotheses of Theorem \ref{unit}. 

In particular, these hypotheses can be realized over a finite field,
since they impose three conditions on the $4$-dimensional
variety of involutions in $PGL_3$.
In fact, over a finite field more is true.

\begin{theorem} (= Theorems \ref{4.18} and \ref{5.18})
Suppose that $k$ is a finite field
and that $g$ is a hyperbolic element of $Cr_2(k)$.
Then $g$ is tight and $g\not\in\langle\langle g^N\rangle\rangle$
for all sufficiently divisible $N$.
\noproof
\end{theorem}

\end{section}
\begin{section}
{Background and motivation: dynamical systems and
symmetric key cryptography}\label{crypto}
H{\'e}non introduced certain complex quadratic plane
Cremona transformations,
now called H{\'e}non maps, as models of (sections of) dynamical systems
such as the Lorenz equations.
They are of the form 
$$f(x,y)=(ay+q(x),x),$$
where $q$ is a quadratic polynomial and $a$ is a non-zero scalar.

Then $f$ might have sensitive dependence on initial conditions 
in this sense:
%
even if initial points $x_0$ and $y_0$ are very close,
their images $f^n(x_0)$ and $f^n(y_0)$ can be far apart
for large values of $n$. 

In the context where $f$ is a smooth self-map of a compact 
manifold $X$ this, when stated precisely in terms
of Lyapunov exponents, turns out to be equivalent to 
the topological entropy $h(f)$ being strictly positive.
(Recall that the topological entropy $h(g)$ of a self-map
$g$ of a compact metric space $X$ is defined by
$$h(g)=\lim_{\eps\to 0}\limsup_{n\to\infty}\log\left(\frac{1}{n}N(n,\eps)\right),$$
where $N(n,\eps)$ is the number of $g$--orbit segments of length $n$ that
are at least a distance $\eps$ apart. Gromov [G] extended
this definition to cover correspondences, which include Cremona
transformations, as well.) 
However, this definition involves two limits, so the questions arise
of finding how large $n$ must be taken, and how small $\eps$,
in order to estimate it in to a given accuracy in a bounded time.


On the other hand,
over a finite field, especially one of characteristic $2$,
Feistel introduced the same kind of Cremona transformations
$$f(x,y)=(y+q(x),x),$$ 
except that he took the parameter
$a$ to be $a=1$ always and he did not demand that $q$ be quadratic.
Note that, in characteristic $2$, this map $f$ is the composite
of two involutions: 
$$f=\a\circ\s,$$
where $\a$ is the linear map $(x,y)\mapsto(y,x)$
and $\s(x,y)=(x,y+q(x))$. 
These maps are also known as \emph{round functions} 
and they are an essential element of \emph{Feistel ciphers} such as
DES, the Data Encryption
Standard. The Advanced Encryption Standard, AES,
uses different Cremona transformations, but otherwise both DES and AES
have a similar structure. In fact, as explained below,
AES is a Cremona transformation that is an element of a Galois twist of the
group of \emph{standard} Cremona transformations
of $\P^{128}_{\F_2}$. This (the untwisted group,
that is) is the group generated by $PGL_{128}(\F_2)$
and the standard non-linear birational
involution $x_i{\mapstochar\dashrightarrow}  x_i^{-1}$.

Here is a toy model of DES ; it is a toy because it omits
the \emph{key schedule}. Taking the key schedule into account,
as is done below,
gives something like a noisy dynamical system, but where
the noise is wholly determined in advance,
as part of the infrastructure.

After Alice and Bob have established a key $K$
(for example, by using some version of public key cryptography
based, say, on elliptic curves), the key determines,
according to a fixed public procedure
that is part of the infrastructure of the algorithm,
a round function $f_K=\a\circ\s_K$ that is a 
Cremona transformation of some projective space
$\P^n_k$ over a finite field $k$. 

Once the key has been established,
encryption of a message $M$ is this: break $M$ into blocks $M_i$,
each of size $n$ (that is, $M_i$ is a $k$-point of $\A^n$)
and then, for a fixed integer $N$ that is also part of
the infrastructure of the algorithm, apply the transformation
$f_K^N$ to the plaintext block $M_i$ and transmit $f_K^N(M_i)$.
Decryption is: apply $(f_K^{-1})^N$ to each block
$f_K^N(M_i)$ that is received. It is a basic
requirement that, given possession of the key $K$, decryption
should be as fast as encryption; this is achieved by constructing
$f_K$ as the product of two involutions, and then
decryption is merely the process of applying the same two
involutions but in the opposite order.

The key schedule amounts to fixing
$k$-linear involutions $L_1,\ldots,L_N$
as part of the infrastructure (so that, in particular, the $L_i$
are independent of the key), defining
$f_i=f_{K,i}=f_K\circ L_i$ and then taking encryption
to be the iterate $f_N\circ\cdots\circ f_1$.
Since each iterand $f_i$ is a product of three involutions,
decryption is merely a matter of reversing the
order of these involutions,
and so is as cheap, in terms of time and memory, as encryption.

AES can be described in similar terms.
First, some Galois twist $\s$ of the standard non-linear birational involution
$$(x_0,...,x_n){\mapstochar\dashrightarrow} (x_0^{-1},...,x_n^{-1})$$
of $\P^n_k$ is given in advance and is public.
Then, after the key $K$ has been established, as above, it
is used to construct a linear transformation
of $\P^n_k$ and we set $f_K=\s\circ L_K$.
Thereafter the algorithm runs as for DES.
(If the process ever encounters a base point, meaning that
it is trying to invert $0$, then it maps $0$ to $0$.
So the iterated Cremona transformation is garbled.
However, the security of the scheme does not reside in this garbling.)
Since $f_K^{-1}=\s\circ L_K^{-1}$, decrypting is then,
as with DES, as fast and cheap as encrypting (especially if $L_K$
is an involution). However, in higher dimensions, 
the inverse of a general
Cremona transformation $\phi$ is of higher degree than $\phi$,
so that inversion of $\phi$ is slower and more expensive 
than the execution of $\phi$.

As with DES, so AES has a key schedule, and the structure is similar.

Encryption should also mix up the points of projective space thoroughly
and quickly; in other words, it is desirable that if $x_0$ and $y_0$
are distinct basepoints that are close, then the points
$f_K^N(x_0)$ and $f_K^N(y_0)$ should be far apart
for some large, but fixed, value of $N$. In other words,
the round function should be sensitive
to the parameters that define it, in the sense of (1) above.
That is, in highly simplified terms,
over the complex numbers
certain Cremona transformations 
serve as simple models of a process that is known to be chaotic,
while over a finite field the same Cremona transformations
are used to create a process that merely
has a convincing appearance of chaos.

Moreover [M1], positive entropy
does not exclude the existence of Siegel discs; Siegel discs
(whatever their analogues might be over finite fields)
are undesirable in a cryptographic context because they are regions
consisting of
plaintexts that are close and that remain close after encryption.
On the other hand, the theorem of Gromov and Yomdin,
that the entropy of an endomorphism $g$ of a smooth projective variety
$X$ is the logarithm of the spectral radius of the action of $g$
on the cohomology of $X$, shows that,
as a consequence of the Lefschetz fixed point formula, $h(g)$ can be computed by
counting fixed points of a certain number of iterates of $g$;
how many iterates are required
depends on the Betti numbers of $X$. 

Despite the fact that the contemporary algorithm
AES is a Cremona transformation,
in this paper 
we consider only plane Cremona transformations.
The reason is simple: we don't know any analogous
results in higher dimensions, even for the 
group of standard transformations.
\end{section}

\begin{section}{Hyperbolic space}
Here we review the construction and basic properties of
the infinite hyperbolic space $\mathfrak H=\mathfrak H_k=\mathfrak H(\P^2_k)$
and the action of $Cr_2(k)$ on it.
This is taken from [CL]; we repeat it only
in order to establish notation.

Let $V$ be any smooth projective surface over the field $k$.
Set $\sZ(V)_\Z =\varinjlim_{Y\to V}\NS(Y)$,
where the direct limit is taken over all blow-ups
$Y\to V$, and $\sZ(V)=\sZ(V)_\Z\otimes_\Z\R$.
There is a hyperbolic completion $\tsZ(V)$ of
$\sZ(V)$, given by
$$\tsZ(V)=\{\lambda +\sum_{e_P\in\sE} n_Pe_P\vert \lambda \in \NS(V)_\R,                              
n_P\in \R, \sum \deg (P)n_P^2<\infty\},$$
where $e_P$ is the exceptional curve associated
to the closed point $P$ on some blow-up $Y$,
$\sE$ is the set of such curves,
$\deg(P)$ is the degree of the field extension $k(P)/k$
and $\sum\deg(P)n_P^2<\infty$
means that, for all $\eps>0$ there is a finite subset $\sF$ of $\sE$
such that, for all finite subsets $\sG$ of $\sE-\sF$, we have
$\sum_{e_P\in\sG}\deg(P)n_P^2<\eps$.
(By definition, $\sum \deg (P)n_P^2$ is the number $l$
such that for all $\eps>0$
there is a finite subset $S$ of closed points $P$
such that, whenever $S\subset T$ and $T$ is a finite set of closed points,
$\big\vert \sum_{P\in T}\deg (P)n_P^2 -l\big\vert<\eps.$)
On $\tsZ(V)$ there is a hyperbolic inner product
denoted by $(x.y)$.
The hyperbolic space $\mathfrak H(V)$ is one of the two connected
components of the locus $\{x\in \tsZ(V)\vert (x.x)=1\}$.
Note that, if $x=\lambda +\sum n_Pe_P$, then
$(x.x)=(\lambda.\lambda)-\sum_P\deg(P)n_P^2$.
The distance on $\mathfrak H(V)$ is denoted by $d$, so that
$\cosh d(x,y)=(x.y)$. The isometries
of $\mathfrak H(V)$ are the continuous linear transformations
of $\tsZ(V)$ that preserve the inner product
and the connected component above.

Given any blow--up $Y\to V$, there are natural isomorphisms
$\sZ(V)_\Z\to \sZ(Y)_\Z$, $\mathfrak H(V)\to \mathfrak H(Y)$,
etc., so that the group
$\Bir(V)$ of birational automorphisms of $V$ acts
as a group of isometries of $\mathfrak H(V)$
via $g\mapsto g_*$.

From now on, we take $V=\P^2_k$
and write $\mathfrak H(V)=\mathfrak H_k$
and $\Bir(V)=Cr_2(k)$. 

\begin{lemma}\label{subfield}
The formation of $\mathfrak H_k$ is functorial in $k$
and, 
if $k$ is a subfield of $K$, then
$\frak H_k$ is naturally a closed geodesic subspace
of $\frak H_K$ and, as a subgroup of $Cr_2(K)$, 
$Cr_2(k)$ acts on $\mathfrak H_K$ so as to preserve
$\mathfrak H_k$.
\begin{proof} Immediate from the construction
of $\mathfrak H_k$.
\end{proof}
\end{lemma}

We shall usually
drop the subscript $k$ from $\mathfrak H_k$.

Isometries $\phi$ of $\mathfrak H$ are of three types:
$\phi$ is \emph{elliptic} if it has a fixed point in
(the interior of) $\mathfrak H$; $\phi$ is \emph{parabolic}
if it has a unique fixed point on the ideal boundary
$\partial\mathfrak H$ of $\mathfrak H$ and no fixed point in $\mathfrak H$;
$\phi$ is \emph{hyperbolic} if the lower bound
$L(\phi)=\liminf d(x,\phi(x))$ is strictly positive
and is attained in $\mathfrak H$.

In this last case the set $\{x\in\mathfrak H\ \vert\  d(x,\phi(x))=L(\phi)\}$
is a geodesic in $\mathfrak H$. It is the \emph{axis} of $\phi$
and is denoted by $\Ax(\phi)$.
It is the unique geodesic preserved by $\phi$ and its endpoints
on $\partial\mathfrak H$ are the unique fixed point on the closure
$\overline{\mathfrak H}=\mathfrak H\cup \partial\mathfrak H$.
The quantity $L(\phi)$ is the \emph{translation length}
of $\phi$. On the other hand,
a parabolic isometry does not preserve any geodesic.

\begin{lemma}\label{subspace}
Suppose that $\mathfrak G$ is a closed hyperbolic subspace of $\mathfrak H$
and $\phi$ an isometry of $\mathfrak H$ that preserves
$\mathfrak G$.

\part[i] $\phi$ is elliptic if and only if
$\phi\vert_{\mathfrak G}$ is elliptic.

\part[ii]
$\phi$ is hyperbolic if and only
if $\phi\vert_{\mathfrak G}$ is hyperbolic,
and in this case $\Ax(\phi\vert_{\mathfrak G})=\Ax(\phi)$.

\part[iii]
$\phi$ is parabolic if and only if
$\phi\vert_{\mathfrak G}$ is parabolic, and in this case
the unique ideal boundary point of $\mathfrak H$
that is fixed by $\phi$ equals the ideal
boundary point of $\mathfrak G$ that is fixed by $\phi\vert_{\mathfrak G}$.
\part[iv]
$L(\phi\vert_{\mathfrak G})=L(\phi)$.
\begin{proof}
\DHrefpart{i}: If $P\in\mathfrak H-\mathfrak G$ and $\phi(P)=P$,
then $\phi$ also fixes the unique point $Q$ on $\mathfrak G$
that is closest to $P$.

\DHrefpart{ii}: Assume $\phi$ to be hyperbolic, with axis $\G$.
There is a map $\G\to\mathfrak G:x\mapsto y$
if $y$ is the closest point to $x$ that lies on $\mathfrak G$.
The image is a geodesic $\Delta$ which is preserved by $\phi$.
Since $\phi$ does not preserve two geodesics,
$\G=\Delta$.

Conversely, if $\phi\vert_{\mathfrak G}$ is hyperbolic,
it preserves a geodesic $\Delta$ in $\mathfrak G$,
so that, again by \DHrefpart{i}, $\phi$ is hyperbolic
and its axis is $\Delta$.

\DHrefpart{iii}: Obvious, from \DHrefpart{i} and \DHrefpart{ii}.

\DHrefpart{iv}: Obvious.
\end{proof}
\end{lemma}

\begin{corollary}\label{indep of field}
If $g\in Cr_2(k)$,
then whether it is elliptic, parabolic or hyperbolic and its translation length 
can be calculated after any making any
extension of $k$.
\noproof
\end{corollary}

\begin{lemma}\label{important} Suppose that $\G$
is a geodesic in $\mathfrak H$
and that $\phi$ is an isometry of $\mathfrak H$
that preserves $\G$. Then the following
statements hold.

\part[i] Either
$\phi$ is elliptic and fixes a point on $\G$
or $\phi$ is hyperbolic.

\part[ii] If $\phi$ is hyperbolic then
$\G$ is the unique geodesic preserved
by $\phi$ and equals the axis of $\phi$.

\part[iii] The translation length of $\phi$ acting on
on $\mathfrak H$ equals the translation
length of its restriction $\phi\vert_\G$ to $\G$.
\begin{proof} Take $\mathfrak G=\G$
in Lemma \ref{subspace}.
\end{proof}
\end{lemma}

%

\begin{lemma}\label{crucial}
Suppose that $\s,\a$ are involutions of 
$\mathfrak H$.

\part[i] $\Fix(\s)$ and $\Fix(\a)$
are hyperbolic subspaces of $\mathfrak H$.

\part[ii] $\s$ preserves each geodesic that
is perpendicular to $\Fix(\s)$,
and the same for $\a$.

\part[iii] $\Fix(\s)$ and $\Fix(\a)$
meet in $\mathfrak H$ if and only if $s\a=g$, say, is elliptic.

\part[iv] $\Fix(\s_*)$ and $\Fix(\a_*)$ are parallel
if and only if they meet in a single ideal boundary
point $P$
if and only if $g$ is parabolic.

\part[v] If $\Fix(\s_*)$ and $\Fix(\a_*)$ are ultraparallel
then there is a unique geodesic $\G$ perpendicular to both
and $g$ is hyperbolic.
This geodesic is preserved by both $\a$ and $\s$
and is the axis of $g$.
\begin{proof}
\DHrefpart{i} is clear: both fixed loci
are projectivizations of linear spaces.

\DHrefpart{ii} is a simple observation.

\DHrefpart{iii}: if $\Fix(\s)$ and $\Fix(\a)$ meet in
an interior point $x$ of $\mathfrak H$ then $g(x)=x$
and $g$ is elliptic.

Conversely, suppose $x$ is an interior point and $g(x)=x$.
If $x\in\Fix(\s)\cup\Fix(\a)$ then it is easy to see that
$x\in\Fix(\s)\cap\Fix(\a)$, so suppose
that this is not the case.
Then there is a unique geodesic segment $l$ from
$x$ to $\a(x)$: this is perpendicular to $\Fix(\a)$
and $\Fix(\a)$ cuts $l$ at its midpoint.
Similarly, there is a unique geodesic segment
$m$ from $\a(x)$ to $\s\a(x)$: this
is perpendicular to $\Fix(\s)$
and $\Fix(\s)$ cuts $m$ at its midpoint.
But $\s\a(x)=x$, so $l=m$ and the
midpoint lies in $\Fix(\a)\cap\Fix(\s)$
and \DHrefpart{iii} is proved.

So we can assume that
$\Fix(\s)\cap\Fix(\a)$ contains
no interior point and that $g$ is not elliptic.

\noindent Case (a): $\Fix(\s)\cap\Fix(\a)$ contains
a boundary point $P$. Then $g(P)=P$.
Suppose that $g$ is hyperbolic, with axis $\Ax(g)=l$,
and suppose that the boundary point $Q$ is the other
endpoint of $l$.
Since $\s g\s=g^{-1}$ and $\Ax(g^{-1})=\Ax(g)$,
$\s$ preserves $\Ax(g)$ and so fixes $Q$.
Similarly $\a(Q)=Q$. Then $\Fix(\s)\cap\Fix(\a)$
contains $l$, which is absurd. So $g$ is parabolic.

\noindent Case (b): $\Fix(\s)\cap\Fix(\a)$
contains no interior nor boundary point.
That is, they are ultraparallel.
Then there is a unique geodesic $l$ perpendicular
to both of them. Because $\a,\s$ are involutions,
each of then preserves $l$, and so $g$ is not elliptic
but preserves a geodesic $l$. Then $g$ is hyperbolic
and $l$ is its axis. This proves
\DHrefpart{iv} and \DHrefpart{v}.
\end{proof}
\end{lemma}

Suppose that $\d$ is a Cremona transformation,
of degree $D$ (in that $\d$ is defined by a net
of homogeneous polynomials of degree $D$).
Then [CL] $L(\d_*)$ equals the dynamical degree
of $\d$, defined as
$\lim_{n\to\infty}(\deg(\d^n))^{1/n}.$
Say that 
$\d$ is elliptic, etc.,
if $\d_*$ is so and that
$\d$ is \emph{biregular}
if there is some rational surface $X$ on which $\d$ is
a biregular automorphism.

Over an algebraically closed field
$\d$ is elliptic
if and only if $\d$ is biregular, and in this case
there is a rational surface $X$ on which $\d$ is biregular
and an integer $n>0$ such that $\d^n$ lies in
the connected component $\Aut_X^0$ of the identity element
in $\Aut_X$. 
The map $\d$ is parabolic if and only if it preserves
a pencil of curves of genus at most $1$.

Over any field,
if $\d$ is hyperbolic then there is no
pencil of curves
that is preserved by $\d$.

Moreover, over $\C$, the entropy $h(\d)$
of $\d$ satisfies $h(\d)\le L(\d_*)$
and equality holds if $\d$ is biregular on some
smooth projective rational surface.
\end{section}
\begin{section}{Graphs and lattices}
Fix homogeneous co-ordinates $x,y,z$ on $V=\P^2_k$.
We denote by $\s$ the standard quadratic
involution $\s:(x,y,z){\mapstochar\dashrightarrow} (yz,xz,xy)$ 
and by $\a$ a linear involution. We put $g=\s\a$.

Say that $P,Q,R$ are the base points of $\s$ and $Y=\Bl_{P,Q,R}V$,
with exceptional curves $e_P,e_Q,e_R$. Note that $\s$ is biregular
on $Y$ and $\a$ is biregular on $V$. Let $l$ denote the class
of a line in $V$. 

Let $E$ denote the set of
all exceptional curves $e_x$ as $x$ runs over all closed
points of all blow-ups of $Y$.

Since $\sZ(Y)_\Z\cong\sZ(V)_\Z$, the
lattice $\sZ(V)_\Z$ has a $\Z$-basis
$\{l\}\cup\{e_P,e_Q,e_R\}\cup E$.

From now on we shall not always be careful to
distinguish between $\a$ and $\a_*$,
nor between $\s$ and $\s_*$.

\begin{lemma} $\a$ permutes the set $\{e_P,e_Q,e_R\}\cup E$
and $\s$ permutes $E$.
\begin{proof} Immediate, from the facts that $\a$ is biregular on
$V$ and that $\s$ is biregular on $Y$.
\end{proof}
\end{lemma}

Let $v_0$ denote the
root $v_0=l-e_P-e_Q-e_R$.

\begin{lemma} $\a$ preserves $l$ and $\s$ acts on the lattice
$\Z.\{l,e_P,e_Q,e_R\}$ as the reflection $s_{v_0}$ in $v_0$.
\noproof
\end{lemma}

Our aim is to construct a bipartite graph $H$ that depends 
upon $\a$. Then $g=\s\a$ will act as something close
to a Coxeter element in the corresponding Coxeter group.

We begin by constructing subsets $\tG_\a$ and $\tG_\s$
of the $\Z$-lattice spanned by $\{l,e_P,e_Q,e_R\}\cup E$, as follows.
These subsets are not necessarily disjoint.

The elements (or vertices) of $\tG_\s$ are $v_0$ and one representative $e_x-\s(e_x)$
of each non-zero pair $\pm(e_x-\s(e_x))$ as $e_x$ runs over $E$.
The vertices of $\tG_\a$ are one representative $e_y-\a(e_y)$
of each non-zero pair $\pm(e_y-\a(e_y))$ as $e_y$ 
runs over $E\cup \{e_P,e_Q,e_R\}$. Finally, if $v$ lies in
$\tG_\a$ and $\pm v$ lies in $\tG_\s$, then choose 
$v$ rather than $-v$ in $\tG_\s$.

Put $\tG=\tG_\a\cup\tG_\s$.
Because $\a$ and $\s$ are involutions,
two different points in
$E\cup\{e_P,e_Q,e_R\}$, (resp., in $E$),
cannot give the same vertex in $\tG_\a$ (resp., in $\tG_\s$).

We join two distinct vertices in $\tG$ by an edge of multiplicity equal to
their intersection number, if that number is non-zero.
If the intersection number is zero, then the corresponding vertices
remain disjoint. So every edge has multiplicity $\pm 1$ or $2$.
(Since there are no loops, that is, since no vertex is 
joined to itself, $-2$ does not occur.)
Define the valency of a vertex $v$ to be the sum of the absolute values
of the multiplicities of the edges meeting $v$.

\begin{lemma} If $v$ lies in the intersection $\tG_\a\cap\tG_\s$
then $v$ is disjoint from all other vertices in $\tG$.
\begin{proof} Suppose that $v=e_x-\s(e_x)$ and $v=e_y-\a(e_y)$.
Note that $e_x,\s(e_x),e_y$ and $\a(e_y$
are all classes of irreducible curves.
We proceed to consider three cases separately.
\begin{enumerate}
\item $v$ meets $e_z-\s(e_z)$. Since $e_z$ and
$\s(e_z$ are also classes of irreducible curves,
either $e_x=e_z$, and then $v=e_z-\s(e_z)$, or $e_x=\s(e_z)$. 
This latter possibility contradicts the construction of $\tG_\s$.

\item $v$ meets $e_t-\a(e_t)$. We reach a similar conclusion.

\item $v$ meets $v_0=l-e_P-e_Q-e_R$.
Then $e_x\in\{e_P,e_Q,e_R\}$ but $\a(e_P)\ne\s(e_P)$.
So this cannot happen.
\end{enumerate}
\end{proof}
\end{lemma}

Now define $G_\a=\tG_\a-(\tG_\a\cap\tG_\s)$
and define $G_\s$ similarly. 
Note that $v_0$ lies in $G_\s$.

\begin{lemma}\label{key lemma} 

\part[i] There are no edges within either $G_\a$ or $G_\s$.

\part[ii] $G_\a\cup G_\s$ is a bipartite graph $G$. 

\part[iii] The 
vertex $v_0$ is of valency at most $3$
and every other
vertex of $G$ is of valency at most $2$.
\begin{proof}
It is enough to notice that in $\tG$ the vertex $v_0$
has valency at most $3$ and that every other vertex has
valency at most $2$, so that deleting $\tG_\a\cap\tG_\s$
amounts to deleting those vertices that are joined to themselves
and to no other vertex. Equivalently, deleting $\tG_\a\cap\tG_\s$
amounts to deleting all double bonds and the corresponding vertices.
\end{proof}
\end{lemma}

Now define $\tH$ to be the connected component of $G$ that
contains $v_0$. 
 
\begin{lemma} $\tH$ is bipartite and is \emph{either} a tree 
$T_{p,q,r}$ consisting of $v_0$ and three arms of lengths $p,q,r\le\infty$
attached to $v_0$ \emph{or} the union $\Delta_{m,r}$ of
a cycle of finite length $m$
together with an arm $v_0,w_{r-1},...,w_1$ 
of length $r\le\infty$ attached to the cycle at $v_0$.
\begin{proof} Immediate.
\end{proof}
\end{lemma}
\begin{remark} When we speak of the \emph{length}
of an arm, we count the vertex to which it is joined.
So, for example, $T_{2,3,5}$ is the $E_8$ diagram,
$T_{p,q,r}$ has a total of $p+q+r-2$ vertices
and, in general, a non-trivial arm has length at least $2$.
So $\Delta_{m,1}$ is a cycle of length $m$.
\end{remark}
\begin{lemma} If $\tH=\Delta_{m,r}$ then $m$ is even.
\begin{proof} 
Any subgraph of a bipartite graph is bipartite,
so the cycle in $\tH$ is bipartite, so even.
\end{proof}
\end{lemma}

Write $m=2n$.
From $\tH$, construct a graph $H$ with the same vertices as $\tH$,
but where every edge except at most one has multiplicity $+1$,
by starting at $v_0$ and proceeding either outwards along
one arm at a time (in the case of $T_{p,q,r}$) or around the cycle
and then along the arm (in the case of $\Delta_{2n,r}$) as follows:
at each step, change an edge of multiplicity $-1$ into an
edge of multiplicity $+1$ by replacing a vector $v$ by its
negative, $-v$. If $H$ is of type $T_{p,q,r}$ then all its
edges have multiplicity $+1$, and we write $H=T_{p,q,r}$;
in the other case all edges, 
except possibly one edge that meets $v_0$ and lies in the cycle,
are of multiplicity $+1$, and we write $H=\Delta_{2n,r}^\pm$
accordingly.

\begin{lemma}\label{diagram} 
\part[i] $H$ is either of type 
$T_{p,q,r}$ with $p,q,r\le\infty$ 
or of type $\Delta_{2n,r}^-$ with 
$2\le n<\infty$ and $2\le r\le\infty$.

\part[ii] If $p,q,r$ are finite and the points
$P_0,...,P_{p-1},Q_0,...,Q_{q-1},R_0,...,R_{r-1}$
are in general position, then $H$ contains the
diagram $T_{p,q,r}$.
\begin{proof} 
\DHrefpart{i}: 
We use the notation of the preceding proof, and in addition let
$s'$ denote the element of $\{\a,\s\}$ distinct from $s$.
Then there are consecutive nodes in the cycle that are of the form
$e_{sw(P)}-e_{w(P)}$,
$e_{s'sw(P)}-e_{sw(P)}$ and
$e_{sw(Q)}-e_{w(Q)}$,
and also $e_{s'sw(P)}-e_{sw(P)}=e_{s'sw(Q)}-e_{sw(Q)}$.
However, this is absurd.

From the definition of $\Delta_{2n,r}^-$
it is clear that $n,r\ge 2$.

\DHrefpart{ii} is an immediate observation.
\end{proof}
\end{lemma}

Denote by $\Lambda(H)$ the lattice on the vertices $v,w$ of $H$, with pairing
given by the usual intersection pairing of curves.

Put $H_\a=H\cap G_\a$ and $H_\s=H\cap G_\s$,
so that $v_0\in H_\s$.
For $v\in H$, let $s_v$ denote the reflection in $v$.

\begin{lemma} $\s$ acts on $\Lambda(H)$ via the product
$S=\prod_{v\in H_\s}s_v$
and $\a$ acts on $\Lambda(H)$ 
via the product
$A=\prod_{w\in H_\a}s_w$.
\begin{proof} Immediate observation.
Notice that because all the reflections in each product
commute with each other, the order in which they
are taken is immaterial. The fact that each product contains
infinitely many factors is also immaterial, since there
are only finitely many terms in either product that
act non-trivially on any given element of $\Lambda(H)$.
\end{proof}
\end{lemma}

The next lemma is well known but for lack
of a convenient reference we give a proof.

\begin{lemma} \label{degen*}
Suppose that $2\le p\le q\le r<\infty$ and
$H=T_{p,q,r}$. 
Then $\Lambda(H)$ is negative definite if
$1/p+1/q+1/r>1$, degenerate if $1/p+1/q+1/r=1$ and
hyperbolic if $1/p+1/q+1/r<1$.
\begin{proof} If $1/p+1/q+1/r>1$ then $H$ is one of the Dynkin diagrams
classified in Bourbaki [GrLie4-6] and $\Lambda(H)$ is
the corresponding root lattice, twisted by $(-1)$.

If $1/p+1/q+1/r=1$ then $H$ is an affine Dynkin
diagram of type $\tE_{n-1}$, where $n=p+q+r-2$, 
and is degenerate. The radical $R$ is of rank $1$
and $\Lambda(H)/R$ is isomorphic to the root lattice $E_{n-1}$.

If $1/p+1/q+1/r<1$ then $H$ contains $H'=T_{p,q,r-1}$.
By induction, $\Lambda(H')$ is either degenerate or
hyperbolic, and then $\Lambda(H)$ is hyperbolic.
\end{proof}
\end{lemma}

\begin{lemma}\label{degen} $\Lambda(H)$ is degenerate only when $H$ is 
of finite type $T_{p,q,r}$ with $1/p+1/q+1/r=1$.
\begin{proof} 
Suppose that
$H=\Delta_{2n,r}^-$, that $r\ge 2$ and that the diagram is
$$
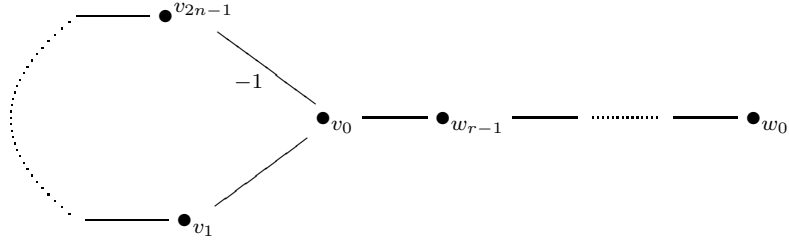
\xymatrix
{
*=0
{}\ar@/_25pt/@{.}[dd]\ar@{-}[r] & {\bullet^{v_{2n-1}}}\ar@{-}[rd]_{-1}\\
&& {\bullet_{v_0}} \ar@{-}[r] & {\bullet_{w_{r-1}}}\ar@{-}[r] &{}\ar@{.}[r] 
& {}\ar@{-}[r] & {\bullet_{w_0}}\\
{}\ar@{-}[r] & {\bullet_{v_1}}\ar@{-}[ru]
}
$$
(the vertices $v_0,\ldots,v_{2n-1}$
are arranged in a cycle of length $2n$). Suppose
that $\eta=\sum_0^{2n-1}a_iv_i+\sum_0^{r-1}b_jw_j$
is in the radical, so that $\eta.v_i=\eta.w_j=0$ for all $i,j$.
Put $w_r=v_0$.

By letting $i$ run from $0$ to $2n-1$, we see that $a_i$
is a linear function of $i$; say $a_i=\lambda i+a_0$
for $i=0,...,2n-1$.
By letting $j$ run from $0$ to $r$,
we see that $b_j=\mu j+b_0$ for $j=0,...,r$.

Since $w_r=v_0$, we get $a_0=\mu r+rb_0$. Also, $b_1=2b_0$, so that
$$\mu=b_0,\ a_0=\mu(r+1).$$

From $\eta.v_{2n-1}=0$ we get 
$a_{2n-2}-2a_{2n-1}-a_0=0,$ so that $2n\lambda +a_0=0$.

From $\eta.v_0=0$ we get
$a_1-a_{2n-1}-2a_0+b_{r-1}=0$, so that
$$0=-2(n-1)\lambda -2a_0+\mu r.$$

Then $\mu=-2\lambda(n+1)/2r$, so that
$\mu=-2n\lambda/(r+1).$ Hence
$(n+1)(r+1)=nr$, which is absurd.

Finally, suppose that $r=1$ (that is,
$H$ is a cycle) and that $\eta=\sum a_iv_i$
is in the radical. Then the equations
$v_i.\eta=0$ for $i=1,...,2n-2$
give 
$$2a_1=a_0+a_2,..., 2a_{2n-2}=a_{2n-3}+a_{2n-1},$$
so that there is some $\lambda$ such that
$a_j=a_0+\lambda j$ for every $j=0,...,2n-1$.
But $v_0.\eta=0$ gives $a_1=2a_0+a_{2n-1}$, 
so that $0=2a_0+(2n-2)\lambda$,
while $v_{2n-1}.\eta=0$ gives
$0=2a_0+2n\lambda$. Then $a_0=\lambda=0$
and $\eta=0$.
\end{proof}
\end{lemma}

There is an obvious natural homomorphism $\beta:\Lambda(H)\to\sZ(V)_\Z$
of lattices.

\begin{lemma} $\beta$ is injective.
\begin{proof}
Inspection.
\end{proof}
\end{lemma}

Let $\Lambda(T_{p,q,r}^{(\lambda)})$ and $\Lambda(\Delta_{2n,r}^{-(\lambda)})$
denote the lattices corresponding to the diagrams
$T_{p,q,r}$ and $\Delta_{2n,r}^-$, but where each vertex
$v$ has $v^2=-\lambda$ and the other intersection numbers
are unchanged. So, for example, 
$\Lambda(H^{(2)})=\Lambda(H)$.

\begin{lemma} 
$\Lambda(\Delta_{2n,r}^{-(2)})$ is hyperbolic if
$2/n+1/r<1$.
\begin{proof} 
$\Lambda(\Delta_{2n,r}^{-(2)})$ is non-degenerate, by Lemma \ref{degen}.
Deleting $v_0$ leaves
a negative definite lattice, and so $\Lambda(\Delta_{2n,r}^{-(2)})$
is either negative definite or hyperbolic.
However, deleting the
vertex opposite $v_0$ in the cycle leaves a $T_{n,n,r}$
diagram, which is hyperbolic.
\end{proof}
\end{lemma}

Say that a lattice $\Lambda$ is \emph{affine} if $\Lambda$ is degenerate,
its radical $R(\Lambda)$ is of rank $1$ and $\Lambda/R(\Lambda)$
is negative definite.
For example, $\Lambda(T_{p,q,r}^{(2)})$ is affine if and only
if $1/p+1/q+1/r=1$.

As $\lambda$ varies, so $\Lambda(T_{p,q,r}^{(\lambda)})$ 
sweeps out a line in the space
of real quadratic forms in $n=p+q+r-2$ variables.
Moreover, $\Lambda(T_{p,q,r}^{(\lambda)})$ is negative
definite if $\lambda \gg 0$
and positive definite if
$\lambda\ll 0$
and the signature of $\Lambda(T_{p,q,r}^{(\lambda)})$
is a non-increasing function of $\lambda$.
(We have adopted the convention that
the signature of a positive definite form
of rank $n$ is $n$ and the signature
of a negative definite form of rank $n$ is $-n$.)
For critical values of $\lambda$ the from 
will be degenerate.

In particular therefore,
there exists $\mu=\mu(p,q,r)$ 
such that $\Lambda(T_{p,q,r}^{(\lambda)})$ is negative
definite if $\lambda>\mu$ while $\Lambda(T_{p,q,r}^{(\mu)})$ is negative
semi-definite and degenerate.
For example, $\mu(2,3,6)=\mu(2,4,4)=\mu(3,3,3)=2$.
Also, define $m=m(p,q,r)$ by 
$$\mu^2=m+m^{-1}+2$$
and $m\ge 1$.

\begin{lemma}\label{increasing}
\part[i] $\mu(p,q,r)$ is a strictly increasing
function of each of $p,q,r$.
That is, $\mu(p,q,r)>\mu(p,q,r-1)$
and the same for $p$ and $q$.
\part[ii] If $1/p+1/q+1/r<1$ then $\mu(p,q,r)>2$.
\begin{proof} \DHrefpart{i}: 
Say $s=\mu(p,q,r-1)$, so that
$\Lambda(T_{p,q,r-1}^{(s)})$ is affine. Then there is a non-zero vector
$\sum n_iv_i$ in the radical of $\Lambda(T_{p,q,r-1}^{(s)})$,
and it is easy to see that we can take every $n_i$ to be
positive. 

Say that $w$ is the vertex adjoined in passing from $T_{p,q,r-1}$
to $T_{p,q,r}$ and that $w$ meets $v_{r-1}$ in $T_{p,q,r-1}$. Then
$$\left(\sum n_iv_i+\epsilon w\right)^2=2n_{r-1}\epsilon -\epsilon^2s>0$$
for $0<\epsilon \ll 1$, so that $\Lambda(T_{p,q,r}^{(s)})$ is hyperbolic
By the discussion above, the act of increasing $s$ will
lead to the lattice becoming negative definite; there is
therefore a critical point $t=\mu(p,q,r)$ at which the lattice becomes
degenerate, and $t>s$. 

\DHrefpart{ii} is a consequence of \DHrefpart{i} and the
the equalities $\mu(2,3,6)=\mu(2,4,4)=\mu(3,3,3)=2$.
\end{proof}
\end{lemma}

Fix $p\le q\le r<\infty$, with $1/p+1/q+1/r<1$.
Say $\mu(p,q,r)=\mu$.

Write $\mu(p,p,p)=\mu_p$ and $m(p,p,p)=m_p$.

\begin{lemma}\label{convergence}\label{3.15}
$\lim_{p,q,r\to\infty}\mu(p,q,r)=3/{\sqrt 2}$
and 
$\mu_p\ge \frac{3}{{\sqrt 2}}(1-2^{-p/4})$.
\begin{proof} Say that $f_1,...,f_{p-1}; g_1,...,g_{q-1}; h_1,...,h_{r-1}$
are the vertices of $T_{p,q,r}$, 
reading outwards along the arms from the
central vertex $v_0$. We can suppose that
$p=q=r$ and then, by symmetry, that
$$\xi=b_p(0)v_0+\sum_{i=1}^{p-1} b_p(i)(f_i+ g_i+ h_i)$$
is in the radical of $\Lambda(T_{p,p,p}^{(\mu_p)})$.

Set $b_p(0)=3$ and $b_p(p)=0$. Then
$0=\xi.v_0=-3\mu_p+3b_p(1)$ and we have a recurrence relation
$$0=\xi.f_n=b_p({n-1})+b_p({n+1})-\mu_p.b_p(n)\ \forall\ n\ge 1.$$
Therefore, by the well known formula for the solution
of such recurrence relations,
there are constants $C_p,D_p$ such that
$$b_p(n)=C_p\t_p^n+D_p\delta_p^n,$$
where $\t_p=(\mu_p+{\sqrt{\mu_p^2-4}})/2$ and
$\delta_p=(\mu_p-{\sqrt{\mu_p^2-4}})/2=\t_p^{-1}\le \t_p$.
Since $b_p(p)=0$, this can be written as
$$b_p(n)=E_p(\t_p^{n-p}-\delta_p^{n-p})$$
for some $E_p$ that is independent of $n$.
Since $b_p(0)=3$ and $b_p(1)=\mu_p$
we get
$$\frac{\mu_p}{3}=\frac{b_p(1)}{b_p(0)}=\frac{\t_p^{2p-1}-\t_p}{\t_p^{2p}-1}.$$
Say $\mu=\lim_{p\to\infty}\mu_p$ and
$\t=\lim_{p\to\infty}\t_p$; then
$$\frac{\mu}{3}=\frac{1}{\t}.$$
This leads to $\mu={\sqrt{3}}/2.$
Since $\mu_p$ is an increasing function of $p$,
we get $\mu_p\to({\sqrt{3}}/2)^-.$

For the speed of convergence, 
note that 
$\mu_p=m_p^{1/2}+m_p^{-1/2}$ and $\mu_p=\t_p+\t_p^{-1}$.
So $m_p=\t_p^2$. The formula
for $\mu_p/3$ above then gives
$$m_p^p=\frac{2m_p-1}{2-m_p}.$$
The same equation is satisfied by $m_p^{-1}$.

\begin{lemma}\label{limit_1}\label{3.16} $m_p\to 2^-$ and
$m_p>2-3.2^{-p/2}$, while
$m_p^{-1}\to (1/2)^+$ and
$m_p^{-1}<1/2 +2^{-p}$.
\begin{proof}
Since $\mu_p\to(3/\sqrt{2})^-$,
$m_p\to 2^-$.

Solving the above equation for $m_p$
when $p=4$ gives $m_4\approx 1.72208,$ so $m_4>\sqrt{2}$.
Since $m_p$ is an increasing function of $p$,
$m_p>\sqrt{2}$ for all $p\ge 4$.

Say $2-m_p=\eta$. Then the equation above gives
$$\eta=(3-2\eta)(2-\eta)^{-p}.$$
Since $m_p>\sqrt{2}$, this leads to
$$\eta<(3-2\eta)2^{-p/2}<3.2^{-p/2}.$$
The proof for $m_p^{-1}$ is similar.
\end{proof}
\end{lemma}
Since $\mu_p^2=m_p+m_p^{-1}+2$, we get
$$\mu_p^2>2-3.2^{-p/2}+\frac{1}{2}+2=3\left(\frac{3}{2}-2^{-p/2}\right).$$
So $\mu_p>\frac{{\sqrt{3}}}{2}(1-\frac{2}{3}2^{-p/2})$
and Lemma \ref{convergence} is proved.
\end{proof}
\end{lemma}

Fix $p\le q\le r<\infty$, with $1/p+1/q+1/r<1$.
Say $\mu(p,q,r)=\mu$.

Then there is a unique totally isotropic vector
$\xi=\xi(p,q,r)\in\Lambda(T_{p,q,r}^{(\mu)})$
where the coefficient of the branch vertex is $1$.

Fix a subset $E$ of $\{p,q,r\}$ and
let the members of $\{p,q,r\}-E$ tend to $\infty$.
Let $H_E$ denote the corresponding infinite graph.

\begin{lemma}\label{new cgnce}
$\mu\to\mu_E$ for some $\mu_E\le 3/{\sqrt{2}}$ and
in the completed hyperbolic vector space
generated by the vector space $\Lambda(H_E)\otimes \R$
there is an eigenvector  $\xi$,
with eigenvalue $\mu_E$,
of the adjacency matrix of $H_E$.
\begin{proof} 
The existence of $\mu_E$
follows from Lemmas \ref{increasing}
and \ref{convergence}.

To prove that $\xi$ exists,
we write it down. 
Define $\t=(\mu_E+{\sqrt{\mu_E^2-4}})/2$.
Suppose, for the sake of explicitness,
that $E=\{p,q\}$. Then
$$\xi=\sum_0^pa_ie_i + \sum_0^qb_jf_j+\sum_0^\infty c_kg_k,$$
where $e_0=f_0=g_0=v_0$,
the branch vertex, $a_p=b_q=0$ and $a_0=b_0=c_0=1$.
and we require that $\xi$ be totally isotropic in
$\Lambda(T_{p,q,\infty}^{(\mu_E)})$.
This is achieved by taking
$$a_i= A(\t^i-\t^{-i})+\t^{-i},\ 
b_j=B(\t^j-\t^{-j})+\t^{-j},\ 
c_k=\t^{-k}$$
where $A=1/(1-\t^{2p})$ and $B=1/(1-\t^{2q})$.
\end{proof}
\end{lemma}

\begin{corollary} If $p,q$ are fixed then
$\lim_{r\to\infty}m(p,q,r)$ exists
and $\le 2$. There is a similar limit if $p$ is fixed
and $q,r\to\infty$.
\noproof
\end{corollary}

Now suppose that $H=\Delta_{2n,r}^-$
and that
$2/n+1/r<1$. 
Define $\mu^{\Delta}(2n,r)$ in the same way
that $\mu(p,q,r)$ is defined for $T_{p,q,r}$.
\begin{lemma}
\part[i] $\mu(n-1,n-1,r)\le\mu^{\Delta}(2n,r)$.
\part[ii] $\mu^{\Delta}(2n,r)$
is a strictly increasing function of both $n$ and $r$.
\part[iii] $\lim_{r\to\infty}\mu^{\Delta}(2n,r)$ exists and $\le 3/{\sqrt{2}}.$
\begin{proof}
\DHrefpart{i}: delete the three vertices $u_1,u_2,u_3$
consisting of the vertex opposite $v_0$ and the two vertices
adjacent to it. This shows that
$\Lambda(T_{n-1,n-1,r})\subset\Lambda(\Delta_{2n,r}^-)$.
Say $\mu(n-1,n-1,r)=\lambda$, and pick
$\xi$ in the radical of $\Lambda(T_{n-1,n-1,r}^{(\lambda)})$.
Then there is a vector $\eta=\xi+\sum_1^3\a u_i$
in $\Lambda(\Delta_{2n,r}^{-(\lambda)})$
with $\eta^2>0$.
So $\mu(n-1,n-1,r)\le\mu^\Delta(2n,r)$.

\DHrefpart{ii} and \DHrefpart{iii}
are proved by the same kind of calculation used in the proof
of Lemma \ref{3.15}.
\end{proof}
\end{lemma}

\begin{corollary} $\mu^\Delta(2p,p)>2-3.2^{-p+1}$.
\begin{proof} Apply Lemma \ref{3.16}.
\end{proof}
\end{corollary}

\begin{lemma} If $n$ is fixed
and $r=\infty$ then
there is an eigenvector 
$\xi$, with eigenvalue
$\lim_{r\to\infty}\mu^\Delta(2n,r)$,
of the infinite adjacency matrix.
\begin{proof}
Put $\mu_E=\mu^\Delta(2n,\infty)$ and 
$\t=(\mu_E+{\sqrt {\mu_E^2-4}})/2$.
As in the proof of Lemma \ref{new cgnce}
the easiest thing is to write $\xi$ down:
$$\xi=\sum_{i=0}^{2n-1}a_iv_i +\sum_{j=0}^\infty b_jw_j$$
where $b_j=\t^{-j}$
$a_i=A_1\t^i+A_2\t^{-i}$ and $A_1,A_2$ are determined by
the equations
$$0=-s+A_1(\t-\t^{2n-1})+A_2(\t^{-1}-\t^{-(2n-1)})+\t^{-1},$$
$$0=A_1(\t^{2n-2}-s\t^{2n-1})+A_2(\t^{-(2n-2)}-s\t^{-(2n-1)})-1.$$
\end{proof}
\end{lemma}

Suppose now that 
$1/p+1/q+1/r<1$
or $2/n+1/r<1$, as appropriate. Then
at this point we have shown, 
whether $H$ is finite or infinite,
the existence of a totally isotropic vector $\xi$
in $\Lambda(H^{(\mu)})$ for some $\mu >2$.
Moreover, $\xi$ generates the radical $R(\Lambda(H^{(\mu)}))$.

Let $M$ denote the adjacency matrix of $H$;
then the matrix $B=-\mu.1+M$
is the Gram matrix of $\Lambda(H^{(\mu)})$
and
$\xi$ is an eigenvector,
with eigenvalue $0$,
of $B$.
The adjacency matrix is, according to the
bipartite decomposition $H=H_\s\sqcup H_\a$, of the form
$$M=\left[
\begin{array}{cc}
0 & {{}^tC}\\
C & 0
\end{array}
\right],$$
where each row and column of $C$ contains at most three
non-zero entries, and each non-zero entry is $1$.
(In the case of $\Delta_{2n,r}^-$ some entries will be $-1$.)

Since $\xi$ generates $R(\Lambda(H^{(\mu)}))$,
the other eigenvalues of $B$
are real and strictly negative. In particular,
$\mu$ is the maximum eigenvalue of $M$,
and is of multiplicity $1$.

Write
$\xi=\left[\begin{array}{c}u\\z\end{array} \right]$.
Then
${}^tCz=\mu u$ and $Cu=\mu z$.

Recall that $g=\s\a$ and notice that, from the bipartite
description of $H$, $\a$ acts as the matrix
$A=\left[\begin{array}{cc}1&0\\C&{-1}\end{array}
\right]$, while $\s$ acts as the matrix
$S=\left[\begin{array}{cc}{-1}&{{}^tC}\\
0&1\end{array} \right]$.
Then 
$$S\left[\begin{array}{c}{u}\\0\end{array} \right]=
\left[\begin{array}{c}{-u}\\0\end{array} \right],\ 
S\left[\begin{array}{c}{0}\\{z}\end{array} \right]=
\left[\begin{array}{c}{\mu u}\\{z}\end{array} \right],$$
$$A\left[\begin{array}{c}{u}\\0\end{array} \right]=
\left[\begin{array}{c}{u}\\{\mu z}\end{array} \right],\ 
A\left[\begin{array}{c}{0}\\{z}\end{array} \right]=
\left[\begin{array}{c}{0}\\{-z}\end{array} \right].$$
Thus $S$ and $A$ preserve the real $2$-plane $\Pi$ 
based
by $\left\{\left[\begin{array}{c}{u}\\0\end{array} \right],
\left[\begin{array}{c}{0}\\{z}\end{array} \right]\right\}$
and act on $\Pi$ with respect to this basis as the matrices
$\left[\begin{array}{cc}{-1}&{\mu}\\
0&1\end{array} \right]$ and
$\left[\begin{array}{cc}1&0\\
{\mu}&{-1}\end{array} \right]$, respectively.

The matrix of $g_*$ is then
$SA=\left[\begin{array}{cc}{\mu^2-1}&{-\mu}\\
{\mu}&{-1}\end{array} \right]=\delta^2,$
where $\delta=\left[\begin{array}{cc}{\mu}&{-1}\\
1&0\end{array} \right].$
Since $\Tr\delta=\mu>2$, it follows that
$\delta$, and so $g_*$, is hyperbolic.

The projectivization of $\Pi$
is a geodesic $\ell$ in the hyperbolic space $\mathfrak H$
and $\ell$ is preserved by $\s,\a$ and $g_*$.
Therefore,
by Lemma \ref{important}, the translation length of $g_*$ acting on
$\mathfrak H$ equals the translation length
of $g_*$ acting on $\ell$.
We shall calculate this from the preceding discussion.
 
Write $\mu^2=m+m^{-1}+2,$ so that
$$\log m =2\cosh^{-1}(\mu/2).$$

The next result contains the first part of Theorem \ref{theorem 1}.

\begin{theorem}\label{theorem 4.19}\label{4.21} 
Suppose that $g$ is in $(p,q,r)$-general position
and that $1/p+1/q+1/r<1$.
Then $g$ is hyperbolic.
If also $4\le p\le q\le r$
then $L(g)\ge \log (2-3.2^{-p/2})$.
\begin{proof}
The hypotheses imply that $H$ is either
$T_{p,q,r}$ with $1/p+1/q+1/r<1$ or
$\Delta_{2n,r}^-$
with $2/n+1/r<1$.
We now use the notation of the
preceding discussion.

The vector 
$x=\left[\begin{array}{c}1\\0\end{array} \right]\in\Pi$
satisfies $(x.x)=1$ and $(x.\delta(x))=\mu/2$,
so that $d(x,\delta(x))=\cosh^{-1}(\mu/2)$.
That is, $L(\delta)=\cosh^{-1}(\mu/2)$,
so that $g_*$ preserves the geodesic $\ell$.
Therefore
$g_*$ is hyperbolic and 
$$L(g)=2\cosh^{-1}(\mu/2)=\log m>0.$$

For the second part, note that,
since $\mu_E$ is an increasing function of each argument
and, in the notation of Lemma \ref{limit_1}
and the formula immediately preceding it,
$\mu_p^2=m_p+m_p^{-1}+2,$
we have
$L(g)\ge\log m_p>\log(2-3.2^{-p/2})$,
by Lemma \ref{limit_1}.
\end{proof}
\end{theorem}

Finally we give the proof of the rest of Theorem \ref{theorem 1}.

\begin{theorem}\label{theorem a} (= Theorem \ref{Trrr})
Assume that $0<\eps<1/3$. 
Then, if $p\ge 12/\eps$ and
the $3p$ points $P_0,...,P_{p-1},...,R_{p-1}$
are in general position, the Cremona transformation 
$g=\s\a$ is hyperbolic and
$$\log 2-\epsilon< L(g)\le \log 2.$$
\begin{proof}
We know that 
$L(g)\ge \log (2-3.2^{-p/2})$,
so that it is enough to verify that
$\log(2-3.2^{-p/2})\ge\log 2-\eps$
if $p>12/\eps$.
This is elementary.
\end{proof}
\end{theorem}

These techniques also give
Cremona transformations with prescribed properties, as follows.

For example, let us calculate $\mu(2,3,7)$.
In terms of the $T_{2,3,7}$ diagram
{\input Dynkin

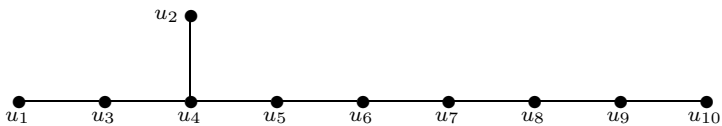
\setdynkin tags=y \setdynkin slope='25
\def\x#1{$\scriptstyle#1$}
\def\present#1{\hbox{\hbox to2cm{$#1$:\hfil}\hskip1.5cm}}
\par\medskip
\dynkinE 10 \x{u_1} \x{u_3} \x{u_4} \x{u_2} \x{u_5} \x{u_6}
\x{u_7} \x{u_8} \x{u_9} \x{u_{10}}
\par
\bigskip
}
\noindent we have to find the value of $s$ for which
the lattice $\Lambda(T_{2,3,7}^{(s)})$ has a non-trivial
totally isotropic vector $\xi=\sum_1^{10}b_iu_i$.
Then $\mu(2,3,7)=s$.
That is, we must eliminate the variables
$b_1,...,b_{10}$ from the equations
$$\xi.u_i=0$$ for all $i=1,...,10$
where $u_i^2=-s$ and the other intersection
numbers are given by the diagram.
We can normalize $\xi$ by assuming that $b_{10}=1$.

This process of elimination is mechanical
and yields
$$s^{10}-9s^8+27s^6-31s^4+12s^2-1=0.$$
Since $\mu(2,3,7)>2$ it follows that
$m(2,3,7)>1$ and
is a zero of Lehmer's polynomial.
Therefore
$$m(2,3,7)=\lambda_{Lehmer}\approx 1.17628.$$
(Recall that $\mu^2=m+m^{-1}+2$.)

Exactly similar calculations show that
$\mu(2,4,5)$ is a solution of
$$s^8-8s^6-20s^4+17s^2-3=0,$$
so that
$m(2,4,5)$ is a zero of
$$m^8-m^5-m^4-m^3+1;$$
$\mu(3,3,4)$ is a zero of
$$s^6-6s^4+8s^2-1$$
so that $m(3,3,4)$ is a zero of
$$m^6-m^4-m^3-m^2+1;$$
$$\mu(4,4,4)={\sqrt{(5+{\sqrt{13}})/2}};$$
$$\mu(5,5,5)={\sqrt{3+{\sqrt 2}}}.$$
The approximate values are given in this table.
\begin{center}
\begin{tabular}{c || c| c| c| c| c| c}
{$(p,q,r)$} & {(2,3,7)} & {(2,4,5)} & {(3,3,4)}&{(4,4,4)} & {(5,5,5)}&{$(\infty,\infty,\infty)$}\\ 
\hline
{$\mu$} & {2.0066} & {2.0153} & {2.0285}&{2.0743}&{2.101}&{$3/\sqrt{2}$}\\ 
\hline
{$m$} & {$\lambda_{Lehmer}$} &{1.28064}&{1.40127}&{1.6644}&{1.8832}&{2}
\end{tabular}
\end{center}

Fix (finite) integers $p,q,r$. Then we can force $H$ to be equal to
the finite tree $T_{p,q,r}$
by putting stringent conditions on the linear involution $\a$, as follows:
\smallskip

\noindent $P_0,...,P_{p-1},...,R_{r-1}$
must be in general position while $P_p=P_{p-1},Q_q=Q_{q-1},R_r=R_{r-1}$.
\smallskip

Note that, since $\s$ has only finitely many fixed points while
$\a$ has a line of fixed points, and (if $\ch k\ne 2$)
one isolated fixed point, 
these are two
conditions on $\a$ (which moves in a $4$-dimensional subvariety
of the $8$-dimensional group $PGL_{3,k}$)
for each of $p,q,r$ that is odd, and
one condition on $\a$ for each that is even. As before,
$l$ is the class of a line in $\P^2$.

\begin{lemma} $H$ is a finite diagram of type $T_{p,q,r}$
and its vertices are $v_0=l-e_{P_0}-e_{Q_0}-e_{R_0};
\ e_{P_0}-e_{P_1},...,e_{P_{p-2}}-e_{P_{p-1}};
\ e_{Q_0}-e_{Q_1},...,e_{Q_{q-2}}-e_{Q_{q-1}};
\ e_{R_0}-e_{R_1},...,e_{R_{r-2}}-e_{R_{r-1}}$.
\begin{proof} Immediate.
\end{proof}
\end{lemma}

\begin{lemma} Both $\a$ and $\s$ are biregular on the blow-up
$Y$ of $\P^2$ at the $p+q+r$ points $P_0,...,R_{r-1}$.
\begin{proof} Also immediate.
\end{proof}
\end{lemma}

Of course, $g$ is then also biregular on $Y$.

\begin{lemma} The orthogonal complement $\Lambda(H)^\perp$ of $\Lambda(H)$
in $\NS(Y)$ has a $\Z$-basis $(x_P,x_Q,x_R)$ given by
$x_P=l-\sum e_{P_i}$, etc. Each of $\a,\s$ acts trivially
on $\Lambda(H)^\perp$.
\noproof
\end{lemma}

\begin{theorem}\label{Coxeter} (= Theorem \ref{theorem 2}) The biregular quadratic map
$g$ acts on $\NS(Y)$ as a Coxeter element in the Weyl
group of the lattice $\Lambda(T_{p,q,r})$.
In particular, if $(p,q,r)=(2,3,5)$ then $g$ has order $30$;
if $(p,q,r)=(2,3,6)$ then $g$ is parabolic;
while if $(p,q,r)=(2,3,r)$ with $r\ge 7$
then $g$ is hyperbolic and
$L(g)=\log\lambda_{r+3}$. In particular,
if $r=7$ then $L(g)=\log\lambda_{Lehmer}$.
\begin{proof} This follows from the discussion above.
\end{proof}
\end{theorem}
\begin{remark} 
Blanc and Cantat prove [BC] that
$L(g)\ge \log\lambda_{Lehmer}$ for any hyperbolic element
$g$ of $Cr_2(k)$, while in [M1] McMullen proves the existence
of biregular Cremona transformations $g$ with 
$L(g)=\log\lambda_{r+3}$ for all $r\ge 7$.
His examples also preserve a cuspidal cubic curve
but are less explicit than ours.
For example, sufficient conditions on the linear involution
$\a$ for $g$ to be biregular of type $(2,3,7)$
are that $\a(P)$ should be a fixed point of $\s$ and that 
$\s\a(Q)$ and $(\s\a)^3(R)$
should lie on the line that is the $1$-dimensional
part of $\Fix_\a$. These are $4$ explicit 
polynomial conditions on the $4$-dimensional
family of involutions in $PGL_3$. 
\end{remark}
\begin{remark}
There is a much shorter argument that suffices
to prove merely that $g$ is hyperbolic
if the diagram $H$ is hyperbolic, without any
estimates, as follows.

Consider the action of the involutions $\s,\a$
on the completed hyperbolic space $\mathfrak H$
associated to $\Lambda(H)$.

\begin{lemma} $\Fix(\a_*)$ and $\Fix(\s_*)$
are ultraparallel.
\begin{proof}
$\Fix(\a_*)\cap\Fix(\s_*)$ is
spanned by common non-zero eigenvectors
$p=\left[\begin{array}{c}x\\y\end{array} \right]$
of $\a_*$ and $\s_*$ such that $(p.p)\ge 0$.
Here, $x$ (resp., $y$) lies in the
Hilbert space completion of the negative definite
vector space $\Lambda(H_\a)\otimes\R$
(resp., $\Lambda(H_\s)\otimes\R$),
each of which is naturally embedded in
the hyperbolic completion of $\Lambda(H)\otimes\R$.

We shall check that no such vectors $p$ exist.
There are three cases to consider.
\begin{enumerate}
\item $\a_*(p)=p$ and $\s_*(p)=p$.
Then $Cx=2y$
and ${}^tCy=2x$.
So $(-2+M)p=0$. Since $-2+M$ is
the Gram matrix of $\Lambda(H)$,
$p$
is then orthogonal to every vector
in $\Lambda(H)$.
However, $\Lambda(H)$ is hyperbolic, so non-degenerate.
\item
$\a_*(p)=-p$.
Then $x=0$ and $y\ne 0$, so that
$(p.p)=-2{}^ty y<0$.

\item $\s_*(p)=-p$.
Then $y=0$ and $x\ne 0$,
so that $(p.p)=-2{}^tx x<0$.
\end{enumerate}
\end{proof}
\end{lemma}
Lemma \ref{crucial} now finishes this shorter argument.
\end{remark}
\end{section}

\begin{section}{Rigidity and tightness}\label{7}
Recall from [CL] the crucial notions of \emph{rigidity}
and \emph{tightness}: given $\eps,B>0$, a
hyperbolic conjugacy class $C$ in $G=Cr_2(k)$ is \emph{$(\eps,B)$-rigid}
if
$$\diam(\Tub_\eps\Ax(g)\cap\Tub_\eps\Ax(h))\le B$$
whenever $g,h\in C$ and $\Ax(g)\ne\Ax(h)$. If $\eps'>\eps$ and
$C$ is $(\eps,B)$-rigid then it is also
$(\eps',B')$-rigid for some explicit
$B'=B'(\eps,\eps',B)$ ([CL], 2.3.2). So we can speak of rigidity
without reference to the precise values of
$\eps$ and $B$.

A hyperbolic element is rigid if its conjugacy class is rigid.

\begin{lemma} If $g$ is rigid then so is $g^n$, for every $n\ne 0$.
\begin{proof} Taking $n$th powers of elements in a conjugacy class
gives a new conjugacy class but does not change the set of axes.
\end{proof}
\end{lemma}

A rigid hyperbolic conjugacy class $C$ is \emph{tight} if whenever
$g,h\in C$ and $\Ax(g)=\Ax(h)$, then $h=g^{\pm 1}$. An element
$g$ is tight if its conjugacy class is tight.

For the rest of this section $g$ will denote a fixed hyperbolic
element of $Cr_2(k)$.


Suppose that, for all $i\in\N$, $\Sigma_i$
is a segment of $\Ax(g)$ of length $i$ and that all the $\Sigma_i$
have the same midpoint.
For $\eps>0$ define
$$V_{\Sigma_i,\eps}                                                                                   
=V_{i,\eps}=\{f\in Cr_2(k)\ \vert\  d(x,f(x))<\eps\ \forall\ x\in \Sigma_i\}.$$

\begin{remark}\label{appeal}
Note that, if $g$ is not rigid, then,
by Prop. 3.3 of [CL], for all bounded
segments $\Sigma$ of $\Ax(g)$ and for all $\eps>0$,
there exists $f\in Cr_2(k)$ such that $d(x,f(x))<\eps$
for all $x\in\Sigma$ while $f$ does not preserve $\Ax(g)$.
So, if $g$ is not rigid, then, for every $i$ and every $\eps>0$
the set $V_{i,\eps}$ contains elements $f$ of $Cr_2(k)$
that do \emph{not} preserve $\Ax(g)$.
In particular, if $g$ is not rigid
then every $V_{i,\eps}$ is infinite.
\end{remark}

We now isolate the main part of the argument and formulate it as
a separate result.

\begin{proposition}\label{above}\label{4.3}
Assume that $k$ is algebraically closed and 
that $V_{i,\eps}$ is infinite for all $i$
and for all $\eps>0$.
There is a positive-dimensional affine algebraic
group variety $S$ over $k$ acting biregularly and effectively on a
$k$-rational surface $Y$ such that,
when $S(k)$ is identified with a subgroup of $Cr_2(k)$,

\part[i] for all sufficiently large $i$
and for all $\eps\ll 1$
$V_{i,\eps}$ is a Zariski dense subset of $S(k)$ and

\part[ii] $g$ normalizes $S(k)$.
\begin{proof}
Fix $i\in I$, $x_0\in\Sigma_{i_0}$ and $\eps>0$.
Put $V_{i,\eps,\le r}=\{f\in V_{i,\eps}\ \vert\ \deg f\le r\}$.
Let $\ell\in\frak H$ be the class of a line in $\P^2$
and suppose that $f\in V_{i,\eps}$. Then
$$d(\ell,f(\ell))\le d(\ell,x_0) +d(x_0,f(x_0)) + d(f(x_0),f(\ell))                                   
\le 2d(\ell,x_0)+\eps.$$
That is, the degree of $f$ is bounded. In other words
there is an integer $D$ such that
$V_{i,\eps}=V_{i,\eps,\le D}$.

Now suppose that $f,h\in V_{i,\eps,\le D}$. Then $fh\in V_{i,2\eps}$,
so that $d(\ell, fh(\ell))\le 2 d(\ell,x_0)+2\eps$.
Since $\cosh^{-1}(\N)$ is discrete in $\R$, it follows that
$\deg(fh)\le D$ if, as we now assume, $0<\eps\ll 1$.
So multiplication gives a map
$$V_{i,\eps,\le D}\times V_{i,\eps,\le D}\to V_{i,2\eps,\le D}.$$

According to [BF], the set of Cremona transformations
whose degree is exactly $d$ is naturally the set of $k$-points of
a reduced quasi-projective
scheme $(Cr_2)_d$, while $(Cr_2(k))_{\le D}=\sqcup_{d\le D}(Cr_2)_d(k)$
has no natural structure as a scheme or algebraic space,
although $(Cr_2(k))_{\le D}$ is naturally a noetherian topological space.
Let $Z_{i,\eps,\le D}$ denote the closure of
$V_{i,\eps,\le D}$ in $(Cr_2(k))_{\le D}$,
so that multiplication defines a map
$$m\ :\ Z_{i,\eps,\le D}\times Z_{i,\eps,\le D}\to Z_{i,2\eps,\le D}.$$
By the noetherian property,
$Z_{i,2\eps,\le D}=Z_{i,\eps,\le D}$ for $0<\eps\ll 1$, so that
$m$ is a map
$$m\ :\ Z_{i,\eps,\le D}\times Z_{i,\eps,\le D}\to Z_{i,\eps,\le D}.$$
We can, and do, suppose that $D$ is minimal with respect to the two properties

\noindent (1) $Z_{i,\eps,\le D}$ is infinite and

\noindent (2) $m(Z_{i,\eps,\le E}\times Z_{i,\eps,\le E})$
is not contained in $Z_{i,\eps,\le E}$ for any $E\le D-1$.

We can write $Z_{i,\eps,\le D}=\sqcup_{d\le D}Z_{i,\eps,d}(k)$,
where the $k$-scheme $Z_{i,\eps,d}$ is the Zariski closure of $V_{i,\eps,d}$
in $(Cr_2)_d$.
\begin{lemma}
$m(Z_{i,\eps, D}(k)\times Z_{i,\eps, D}(k))$
meets $Z_{i,\eps, D}(k)$ non-trivially.
\begin{proof}
Assume otherwise, so that
$m(Z_{i,\eps, D}(k)\times Z_{i,\eps, D}(k))$
lies in $Z_{i,\eps,\le D-1}$. We use
the notation of [BF]: $W_d$ is the projectivized space of triples of
homogeneous degree $d$ polynomials in three indeterminates,
$H_d$ is the locally closed subscheme of $W_d$ consisting
of triples that define a Cremona transformation
of degree at most $d$, and, for any $e\le d$,
$H_{e,d}$ is the locally closed subscheme of $H_d$
of triples that define a Cremona transformation of degree exactly $e$.
There are surjections $\pi_{e,d}:H_{e,d}\to (Cr_2)_e$
and $\pi_d:H_d(k)\to (Cr_2(k))_{\le d}$.
The morphism $\pi_{d,d}$ is an isomorphism.

Set $Z'_{i,\eps, D}=\pi_{D,D}^{-1}(Z_{i,\eps, D})$,
so that there is a commutative diagram
$$                                                                                                    
\xymatrix{                                                                                            
{Z'_{i,\eps, D}(k)}\ar@{^{(}->}[r]\ar[d]^{\cong}&                                                     
{H_{D,D}(k)}\ar@{^{(}->}[r]^{open}\ar[d]^{\cong}&                                                     
{H_D(k)}\ar[d]^{\pi_D}\\                                                                              
{Z_{i,\eps, D}(k)}\ar@{^{(}->}[r]&{(Cr_2)_D(k)}\ar@{^{(}->}[r]&{(Cr_2(k))_{\le D}}                    
}                                                                                                     
$$
From our assumption, there exists $E\le D-1$ and a non-trivial open
subscheme $U$ of $Z_{i,\eps, D}\times Z_{i,\eps, D}$
such that $m(U(k))\subset Z_{i,\eps,E}(k)$. Write
$U'=\pi_{D,D}^{-1}(U)$ and let
$m_D:H_D\times H_D\to H_{D^2}$ denote the multiplication.
Then $m_D(U')\subset H_{E,D^2}$.

Let $\overline{Z'_{i,\eps, D}}$ be the closure of $Z'_{i,\eps, D}$
in $H_D$.
Since $\overline{Z'_{i,\eps, D}}$ maps onto $Z_{i,\eps, D}$ under
$\pi_D$, it follows that
$m(Z_{i,\eps, \le D}\times Z_{i,\eps, \le D})\subset Z_{i,\eps,\le E}$.
In particular,
$m(Z_{i,\eps, \le E}\times Z_{i,\eps, \le E})\subset Z_{i,\eps,\le E}$,
so that, by the minimality assumption,
$Z_{i,\eps, \le E}$ is finite. But then
$m(Z_{i,\eps, \le D}\times Z_{i,\eps, \le D})$ is finite,
which contradicts the fact that
$Z_{i,\eps, \le D}$
is infinite.
\end{proof}
\end{lemma}

Now define $Z_i$ as the closure of $Z'_{i,\eps,D}\cup\{1\}$
in $H_D$. By what we have proved so far, $m$ then defines a rational map
${\overline{m}}:Z_i\times Z_i--\to Z_i$. Since $Z_i$ has an identity element
and is preserved under taking inverses and since multiplication
is associative in $Cr_2$, ${\overline{m}}$ is a group chunk.
By the theorem of Weil and Rosenlicht
there is then a rational surface $Y$ over $k$, a
positive-dimensional group variety
$S$ over $k$ that is $k$-birational to $Z$ and an embedding
of $S$ into the group scheme $\Aut_Y$.
Since $Y$ is rational, the group variety $S$ is affine (it has
no abelian part), and \DHrefpart{i} of Proposition \ref{above} is proved.

For \DHrefpart{ii}, note first that,
if we fix
$i_0\gg 0$ and take any $j,i$ with $j\ge i+L(g)$
and $i\ge i_0$,
then
$V_{i,\eps}\supset V_{j,\eps}$ and $V_{j,\eps}$ is Zariski
dense in $S(k)$.

Let $f\in V_{j,\eps}$; then $f\in V_{i,\eps}$ and
then $gfg^{-1}\in V_{i,\eps}$, since $L(g)\le j-i$.
Therefore there is a
Zariski dense subset
$\Sigma$ of $S(k)$ such that $g\Sigma g^{-1}\subset S(k)$.
It follows that
$g$ normalizes the subgroup $S(k)$ of $Cr_2(k)$.
This concludes the proof of Proposition \ref{4.3}.
\end{proof}
\end{proposition}

\begin{lemma}\label{open}
Assume that $k$ is algebraically closed
and that every $V_{i,\eps}$ is infinite.

\part[i] Every positive-dimensional
subgroup $B$ of $S$ that is normalized by $g$ has an open
orbit in $Y$.
\part[ii] Every positive-dimensional
characteristic subgroup of $S$
has an open orbit in $Y$.
\begin{proof}
\DHrefpart{i}: Suppose that $B$ has no open orbit in $Y$.
Then its orbits form a pencil of curves.
This pencil is then preserved by $g$,
which is impossible since $g$ is hyperbolic.

\DHrefpart{ii} follows at once.
\end{proof}
\end{lemma}

\begin{lemma} (Assume that $k$ is 
algebraically closed and that every $V_{i,\eps}$ is infinite.)

$S$ is not semi-simple.
\begin{proof}
There is an $S$-equivariant desingularization $\tY\to Y$.
Then there is an $S$-equivariant blowing-down
$\tY\to Y'$ where $Y'$ is a minimal rational surface:
that is, either $\P^2$
or a Hirzebruch surface $\Sigma_n$ (that is,
a $\P^1$-bundle $\Sigma_n\to\P^1$
where $n\ge 0$, $n\ne 1$ and $\Sigma_n$
has a unique negative section $\s$ with $\s^2=-n$).
The semi-simple part of the automorphism group
of $\Sigma_n$ is isogenous to $GL_2$
if $n\ge 1$ and is
$PGL_2\times PGL_2$ if $n=0$,
while if $Y'=\P^2$ then $S\subset PGL_3$
and so equals $PGL_2$ if it stabilizes
a conic or equals $PGL_3$ otherwise.

If $S=PGL_2$ then $Y'$ is either $\P^2$
or $\Sigma_n$. If $Y=\Sigma_n$
then the set of rulings on $Y'$ that are preserved
by $S$ is finite and non-empty. Some power of $g$
then fixes at least one of them, which is impossible
since $g$ is hyperbolic.
If $Y'=\P^2$ then $S=PSO_3$ and preserves a conic
$\G$ such that $\G$ and $\P^2-\G$ are the only
$S$-orbits on $\P^2$.
Then every $h\in S$
preserves the base locus of $g$,
so that $g\in PGL_3$.

This argument also covers the case
where $S=PGL_3$.

If $S=PGL_2\times PGL_2$
then $Y'=\P^1\times\P^1$ and then,
by a similar argument,
$g^2$ preserves the two rulings on $Y'$.
But $g$ is hyperbolic.
\end{proof}
\end{lemma}

\begin{lemma} If a torus $T$ acts effectively on a
rational surface $Y$ then $\dim T\le 2$.
\begin{proof} As in the proof of the previous lemma,
we reduce to an action of $T$ on $\P^2$ or $\Sigma_n$.
Then $T$ lies in $PGL_3$ or $PGL_2$
or $PGL_2\times PGL_2$ or a group
isogenous to $GL_2$.
\end{proof}
\end{lemma}

\begin{lemma} 
(Assume that $k$ is algebraically closed
and that every $V_{i,\eps}$ is infinite.)

In $Cr_2(k)$ the element
$g$ normalizes a commutative group $A$
which is isomorphic to
either $\GG_m^2$ or $\GG_a^s$ with $s\ge 2$.
\begin{proof}
Take $H$ to be the connected
component of the centre of the radical of $S$;
since $H$ acts effectively on $Y$ and has an open orbit on $Y$,
it follows that $\dim H\ge 2$.
Therefore $H=\GG_m^r\times\GG_a^s$ with $r+s\ge 2$
and $r\le 2$.

The factor $\GG_a^s$
is characteristic in $S$, so $s\ne 1$,
since each characteristic subgroup has an open orbit in $Y$.
\end{proof}
\end{lemma}

\begin{theorem}\label{lots} $g$ is rigid.
\begin{proof} 
Suppose not, and
assume first that $k$ is algebraically closed.
So every $V_{i,\eps}$ is infinite.
Note that each set $V_{i,\eps}$ contains a subset
$V_{i,\eps,0}$
that is Zariski dense in $S(k)$,
where $S$ and $Y$ are the objects provided by Proposition \ref{4.3}.

Regard $Y$ as the rational surface on which $Cr_2(k)$
acts. In the inverse system of all blow-ups $X\to Y$
there is an inverse subsystem of $S$-equivariant blow-ups;
the normalizer $N$ of $S(k)$ in $Cr_2(k)$
acts on this subsystem. Note that $g$ lies in $N$.

Taking the direct limit of the N{\'e}ron--Severi groups
of the surfaces in this system
and then constructing infinite-dimensional hyperbolic
spaces gives a closed hyperbolic subspace $\frak H_N$
of $\frak H$ that is preserved by $N$. Now $g$
acts hyperbolically on $\frak H$
and therefore acts hyperbolically on every
$g$-invariant closed hyperbolic subspace of $\frak H$.
Therefore $\Ax(g)\subset \frak H_N$.
However, $A(k)$ acts trivially on the
N{\'e}ron--Severi group of every $A$-surface
in the subsystem above, so acts trivially
on $\frak H_N$ and so on $\Ax(g)$. Then
$S(k)$, and so $V_{i,\eps,0}$,
preserves $\Ax(g)$, which,
as in the Remark preceding Proposition \ref{4.3},
contradicts our assumption.

So we have proved that $g$ is rigid when $k$
is algebraically closed.

Now suppose that $k$ is an arbitrary field
and that $g$ is not rigid.

Fix an algebraic closure
$K$ of $k$. Then the hyperbolic space $\frak H_k$ is,
from its definition, a closed
geodesic subspace of $\frak H_K$ that is preserved by $g$.
By assumption, $g$ acts hyperbolically on $\frak H_k$,
so preserves $\Ax(g)$, which is a geodesic $\gamma$ in $\frak H_k$
and so a geodesic in $\frak H_K$. Suppose now that
$g$ is not hyperbolic on $\frak H_K$; then $g$ is not parabolic,
since it preserves $\gamma$, so is elliptic. Suppose that
$P\in\frak H_K$ is a fixed point; then $g$ preserves the hyperbolic plane
$\Pi$ spanned by $\gamma$ and $P$. However, no non-trivial
isometry of $\Pi$ can both fix a point and preserve a geodesic.
So $g$ is hyperbolic on $\frak H_K$ and preserves $\gamma$.
It follows that $\gamma$ is the axis of $g$ when $g$ is regarded as an
isometry of $\frak H_K$.

Pick $\eps,B>0$ such that $g$ is $(\eps,B)$-rigid as an element of
$Cr_2(K)$. Suppose that $f,h$ lie in the conjugacy class
of $g$ in $Cr_2(k)$. Then
$$\diam\left(\Tub_\eps\Ax(f)\cap\Tub_\eps\Ax(h)\right)<B$$
where the tubular neighbourhoods are taken in $\frak H_K$,
and then the same inequality holds in the subspace $\frak H_k$.
That is, $g$ is rigid.
\end{proof}
\end{theorem}

\begin{theorem}\label{4.7}
If $k$ is algebraically closed and
no power of $g$ is tight then
$g$ normalizes a copy of either $\GG_m^2$ or $\GG_a^s$.
\begin{proof} We know that $g$ is rigid.
Fix an orientation
of $\Ax(g)$ and consider the subgroup
$$N=\{f\in Cr_2(k)\vert f(\Ax(g))=\Ax(g)\}$$ of $Cr_2(k)$
and its subgroup $N_+$ of index at most $2$
consisting of transformations that preserve the orientation.
Then there is a group homomorphism $\pi:N_+\to(\R,+)$,
where $h\in N_+$ acts on $\Ax(g)$ as a shift by $\pi(h)$.
Note that an element $h$ of $N_+$ is hyperbolic if and only
if $\pi(h)\ne 0$, and then $\vert\pi(h)\vert=L(h)$.
Set $H=\ker\pi$, so that $H=\{f\in Cr_2(k)\vert f(x)=x\ \forall\ x\in\Ax(g)\}$.

Since [BC] the spectrum of $Cr_2(k)$ (the set of lengths
of its hyperbolic elements) is bounded away from zero, the image
$\im\pi$ of $\pi$ is discrete,
so infinite cyclic. Choose $g_0\in N_+$ such that $\pi(g_0)$ generates
$\im\pi$, so that $N_+$ is a semi-direct product
$N_+=H\rtimes \langle g_0\rangle$.

If $H$ is not finite then every $V_{i,\eps}$ is infinite,
so that, by Proposition \ref{4.3} and results following,
there is a $k$-group variety $S$ normalized by $g$ such that
$S(k)\subset H$ and $S$ contains a subgroup
isomorphic to either $\GG_m^2$ or $\GG_a^2$
that is also normalized by $g$.
So we may, and do, assume that $H$ is finite.

We can write $g=hg_0^m$ for some $h\in H$, $m\in\Z$.
Consider the conjugation action of $g_0$ on the finite group $H$;
then $g_0^s$ centralizes $H$ for some $s>0$. That is,
$g_0^s$ lies in the centre $Z(N_+)$. Also,
$g^s=c.g_0^{ms}$ for some $c\in H$, so that
$g^{st}=c^tg_0^{mst}$ for all $t$; choosing
$t$ to be divisible by the order of $H$ leads to
$g^{st}\in Z(N_+)$. Then $fg^{st}f^{-1}=g^{st}$
whenever $f\in N_+$.

If $f\in N\setminus N_+$, put $g_1=fgf^{-1}$. Then
$g_1=h_1g_0^{-m}$ and $g_1^s=c_1g_0^{-m}$
for some $h_1,c_1\in H$. It follows that
$g_1^{st}=c_1^tg_0^{-mst}$; choosing $t$ as before
gives $g_1^{st}=g_0^{-mst}=g^{-st}$
and $fg^{st}f^{-1}=g^{-st}$, so that
$g^{st}$ is tight.
\end{proof}
\end{theorem}

\begin{proposition}
Suppose that $k$ is algebraically closed
and that $g$ normalizes a copy $A$ of $\GG_a^s$
in $Cr_2(k)$.
Then $s=2$ and $A$
acts biregularly
on $\A^2$ via the standard additive action.
\begin{proof}
We have already observed that $s\ge 2$.

As before, we can find a smooth minimal projective surface $Y$
on which $A$ acts effectively
and biregularly and $Y$ is either $\P^2$ or $\Sigma_n$
with $n\ne 1$.
If a group acts
with a dense orbit
then the stabilizers
are conjugate, so that,
if the group is commutative,
its dimension is $2$.
Therefore, if $s\ge 3$, then
the orbits of $A$ are $1$-dimensional
and preserved by $g$.
However, $g$ is hyperbolic.
So $s=2$.

It remains to show that,
whether $Y=\Sigma_n$ or $\P^2$, there
is an open subvariety, homogeneous under $A$,
which is isomorphic to $\A^2$
on which $A$ acts additively.

If $Y=\P^2$ then $A\subset PGL_3$;
conjugating $A$ into the subgroup represented
by strictly upper triangular matrices shows
that $A$ acts additively on $\A^2$.

If $Y=\Sigma_0=\P^1\times\P^1$
then conjugating $A$ into a suitable
subgroup of $PGL_2\times PGL_2$
shows the same thing.

Recall that, if $n\ge 1$, then 
the automorphism group of $\Sigma_n$
is isogenous to a semi-direct product $L_n\rtimes GL_2$
where $L_n=H^0(\P^1,\sO(n))$ and $A$ is conjugate
to a subgroup of $L_n\rtimes U$
where $U\subset GL_2$ is the group
of strictly upper triangular matrices.
In terms of suitable inhomogeneous
co-ordinates $(x,y)$ on $\Sigma_n$,
$L_n$ is the space of polynomials in $x$
of degree $\le n$, and
the element $(v,0)$ of $L_n\times\{0\}$ acts via
$(v,0)(x,y)=(x,y+v(x))$
while the element $(0,\a)$ of $\{0\}\times U\cong\GG_a$ acts
via $(0,\a)(x,y)=(x+\a,y)$.
Since $A$ has a dense orbit in $Y$,
$A$ is not contained in $L_n$.
It is clear that $(v,0)$ and $(0,\a)$
commute for all $\a\in\GG_a$
if and only if $v$ is a constant polynomial,
and the result is proved.
\end{proof}
\end{proposition}

\begin{proposition} Suppose that
$k=\overline{\F}_p$ and that
$g$ normalizes a copy $A$ of $\GG_a^2$ in $Cr_2(k)$.
Then 
$L(g)$ is an integral multiple of $\log p$.
\begin{proof}
We know that $A$ acts additively on $\A^2$.
That is, $a(x)=a+x$.
Since $g$ normalizes this action,
$g$ is biregular on $\A^2$
and we can write
$gag^{-1}=h(a)$ for $a\in A$.
The equation $h(a)(g(x))=g(a+x)$
can be re-written as
$g(a+x)=h(a)+g(x)$.
Pick an origin $0\in\A^2$;
this provides an isomorphism
$A\to \A^2$.

Setting $x=0$ gives
$g(a)=h(a)+g(0)$.
Set $g^n(0)=r_n$;
then
$$g^n(a)=h^n(a)+r_n.$$
Since $r_n$ is constant,
it follows that $g^n$ and $h^n$ are 
defined by polynomials of the same degrees,
so that $\lambda(g)=\lambda(h)$.
(Recall that $L(g)=\log\lambda(g)$ and
$\lambda(g)=\lim_{n\to\infty}\left(\deg(g^n)^{1/n}\right)$.)
Therefore
it is enough to prove
the proposition for automorphisms
$g$ of the group variety $\GG_a^2$.

Choose a finite subfield $\F=\F_q=\F_{p^n}$ of $k$
over which $A$ and $g$ are defined.
Take the non-commutative polynomial ring
$\F[V]$, the quotient of
$W(\F)[F,V]$ (the Dieudonn{\'e} ring of $\F$)
by the ideal $(F)$.
Then, after fixing an identification
$A=\GG_{a,\F}^{2}$, every $\F$-endomorphism
of $A$ is a $\F[V]$-linear endomorphism $\Phi_1$
of a free $\F[V]$-module $M$ of rank $2$.
Moreover, $\deg(g)=p^{\deg (\Phi_1)}$,
where $\deg (\Phi_1)$ is the maximum of the degrees, with respect to $V$,
of the entries of a matrix that represents $\Phi_1$; this is independent of any
choice of basis.

Let $R$ denote the restriction of scalars
$R=R_{\F/\F_p}A$ and let $r:R\to R$ be the automorphism induced
by $g$. Then $\deg r^m=\deg g^m$ for all $m$ and, after fixing a $\F_p$-isomorphism
$R\to\GG_{a,\F_p}^{2n}$, we can identify $r$ with a $2n\times 2n$ matrix
$\Phi$ over the commutative polynomial ring $\F_p[V]$.

Taking degrees of polynomials, with respect to $V$, defines a discrete
valuation $\deg:\F_p(V)\to\Z$.
We define a norm $v$ on $\F_p(V)$
by $v(f)=p^{-\deg(f)}$.
Let $K=\F_p((V))$ be the $v$-adic completion of $\F_p(V)$,
${\overline{K}}$ an algebraic closure of $K$
and $\C_{\overline{K}}$ the $v$-completion of ${\overline{K}}$;
then $\C_{\overline{K}}$ is algebraically closed, and $v$
extends to a norm on it.

We can regard matrices over $\F_p[V]$ as having entries lying in
$\C_{\overline{K}}$. Then $v$ defines a norm on such matrices
$\Psi,\Omega$ with the properties that
$v(\Psi\Omega)\le v(\Psi)v(\Omega)$
and $v(\Psi+\Omega)\le \max\{v(\Psi),\ v(\Omega)\}$.

Then, by definition and by Gel'fand's theorem on the spectral radius,
$$\lambda(g)=\lambda(r)=\lim_{m\to\infty}(v(\Phi^m)^{1/m})=\sup v(\alpha),$$
where $\alpha\in\C_{\overline{K}}$ runs over the spectrum of $\Phi$.
Since $\Phi$ has finite $\F_p[V]$-rank, this supremum is achieved
by an eigenvalue $\alpha$ that is algebraic over $\F_p(V)$.
Then $\deg(\alpha)\in\Q$, so that
$\lambda(g)=p^{\deg(\alpha)}$ is a rational power $p^{a/b}\ge 1$ of $p$.
On the other hand, by [DF] (see also [BC], Th. 1.2)
$\lambda(g)$ is either a Salem number
or a Pisot number, so is a positive integral power of $p$.
\end{proof}
\end{proposition}

The case where $\ch k=0$, $k$ is algebraically closed and $g$
normalizes a copy of $\GG_a^2$ 
does not arise. For 
then $g$ would be linear,
so not hyperbolic.

\begin{proposition} If $g$ normalizes a
$2$-torus, then its length $L(g)$
equals its length as an isometry of the upper half-plane
and so is the logarithm of a real quadratic unit.
\begin{proof}
The same argument as in the additive case above
shows that it is enough to prove this
for the action of $g$ on the torus itself.
This appears in [BC], p.4.
\end{proof}
\end{proposition}

The next result is now immediate.
As we have already observed, given
any overfield $K$ of $k$,
$L(g)$ is the same, whether taken in $Cr_2(k)$ or in $Cr_2(K)$.

\begin{theorem}\label{4.13} 
Assume that $\ch k=0$ or that
$k$ is algebraic.

Suppose that $L(g)$
is not the logarithm of a quadratic unit;
if $\ch k=p$ suppose also that $L(g)$ is not an
integral multiple of $\log p$.

\part[i] Some power of $g$ is tight.

\part[ii] For all sufficiently divisible $n$,
the normal subgroup $\langle\langle g^n\rangle\rangle$
of $Cr_2(k)$ is proper.
\begin{proof} 
\DHrefpart{i}: The previous results can be applied to show
tightness over $\bark$. 
But enlarging the ground field
is irrelevant, and the result follows.

\DHrefpart{ii} now follows from
Theorem 2.10 of [CL].
(Note the typo in \emph{loc. cit.}:
the phrase ``either $h$ is a conjugate of $g$, or...''
should read ``either $h$ is a conjugate of $g^{\pm 1}$, or...''.)
\end{proof}
\end{theorem}

\begin{corollary}
If $g$ is a hyperbolic quadratic
Cremona transformation and $L(g)\ne\log(\frac{1+{\sqrt{5}}}{2})$
(if $\ch k=2$ then assume also that $L(g)\ne \log 2$)
then some power of $g$ is tight and, for all sufficiently divisible $n$,
the normal closure
$\langle\langle g^n\rangle\rangle$ does not contain $g$ and so
is a non-trivial normal subgroup of $Cr_2(k)$.
\begin{proof} It is only necessary to check that $\frac{1+{\sqrt{5}}}{2}$
is the unique quadratic unit between $1$ and $2$.
\end{proof}
\end{corollary}

{\bf{Now suppose that the ground field $k$ is finite.}}

\begin{proposition} For every point $z\in\frak H$ and every
$r>0$, the set of $f\in Cr_2(k)$ such that
$d(z,f(z))<r$ is finite. In particular, $Cr_2(k)$
acts properly discontinuously on $\frak H$.
\begin{proof} Fix $z,r$ and suppose that
$d(z,f(z))<r$. Let $\ell$ be the class of a line in $\P^2$. Then
$$d(\ell,f(\ell))\le d(\ell,z) +d(z,f(z))+d(f(z),f(\ell)),$$
so that $d(\ell,f(\ell))< 2 d(\ell,z) +r.$
Since $\cosh d(\ell,f(\ell))=(\ell.f(\ell))=\deg (f)$, it follows that
$\deg(f)$ is bounded. Because $k$ is finite,
we are done.
\end{proof}
\end{proposition}

We know that, over any field,
there exist hyperbolic elements $g$ of $Cr_2(k)$; for example,
a hyperbolic element of $SL_2(\Z)$
acting on $\GG_m^2$ will do. Fix a hyperbolic element $g$.

Let $M=\Fix(\Ax(g))$ denote the set of elements
$f$ of $Cr_2(k)$ that fix every point on $\Ax(g)$.

\begin{lemma} $M$ is finite.
\begin{proof}
The argument involving the triangle inequality that was used
in the proof of Proposition \ref{above} shows that $(\ell.f(\ell))$
is bounded independently of $f\in M$, and now finiteness follows again.
\end{proof}
\end{lemma}

\begin{theorem}\label{4.18} Some power of $g$ is tight.
\begin{proof} We know that all powers of $g$ are rigid.
Note that $g$ acts on $M$ by conjugation; choose
$n>0$ such that the conjugation action of $g^n$ is trivial.

Put $N=\Stab(\Ax(g))$; then $N$ is a semi-direct
product $N=M\rtimes\langle\gamma\rangle$,
where $\gamma$ is a hyperbolic element of $N$
whose length is minimal, and $g^n$ is a central
element of $N$. So, for every $h\in N$,
we have $hg^nh^{-1}=g^n$, so that
$g^n$ is tight.
\end{proof}
\end{theorem}

\begin{theorem}\label{5.18}
If $k$ is finite then each sufficiently divisible power of $g$
generates a normal subgroup
of $Cr_2(k)$ that does not contain $g$.
\begin{proof} Apply Theorem 2.10 from [CL] once more.
\end{proof}
\end{theorem}
\end{section}

I am grateful to Anthony Manning, Geoff Robinson, 
Caroline Series and Colin Sparrow
for several valuable conversations and emails, and particularly
to Serge Cantat for pointing out and correcting a significant error.

\bibliography{alggeom,ekedahl}
\bibliographystyle{pretex}
\end{document}

%% file: Dynkin.tex

\catcode`\@=11
\let\e@\expandafter
\def\d@nkinsize{1cm}
\newcount\d@nkincount
\d@nkincount"00
\newif\if@ne

\def\LINE{\hbox to \hsize}
\def\setdynkin#1=#2 {\csname d@nkin@#1\endcsname{#2}}
\def\d@nkin@size#1{\def\d@nkinsize{#1}}
\def\d@nkin@slope#1{\d@nkincount#1}
\def\d@nkin@lines#1{\ifnum#1=1\@netrue\else\@nefalse\fi}
\newif\ifd@dotchoice
\newif\ifd@tags
\def\d@nkin@dotchoice#1{\if#1y\d@dotchoicetrue\else\d@dotchoicefalse\fi
\d@setdots}
\def\d@nkin@tags#1{\if#1y\d@tagstrue\else\d@tagsfalse\fi\d@setdots}
\def\d@setdots{%
\ifd@dotchoice
\def\d@nkinball{\e@\dynkinball\e@{\the\t@@@}}%
\ifd@tags
 \let\d@loopargs\d@loopargsa
 \let\count@rg\count@rga
 \let\init@rg\init@rga
 \let\pr@pare\pr@parea
\else
 \let\d@loopargs\d@loopargsbc
 \let\count@rg\count@rgbc
 \let\init@rg\init@rgbc
 \let\pr@pare\pr@pareb
\fi
\else
\def\d@nkinball{\dynkinball}%
\ifd@tags
 \let\d@loopargs\d@loopargsbc
 \let\count@rg\count@rgbc
 \let\init@rg\init@rgbc
 \let\pr@pare\pr@parec
\else
 \let\d@loopargs\d@loopargsd
 \let\count@rg\count@rgd
 \let\init@rg\init@rgd
 \let\pr@pare\pr@pared
\fi
\fi}

\def\d@nkin@ball#1{\e@\let\e@\dynkinball\csname d@nkin#1\endcsname}

\def\spl@ttoksa#1\@t#2\@t#3\t@end{\t@@@={#1}\t@@@@={#2}\global\t@@={#3}}%
\def\spl@ttoksb#1\@t#2\t@end{\t@@@={#1}\t@@@@={}\global\t@@={#2}}%
\def\spl@ttoksc#1\@t#2\t@end{\t@@@={}\t@@@@={#1}\global\t@@={#2}}%

\def\pr@parea{\e@\spl@ttoksa\the\t@@\t@end}%
\def\pr@pareb{\e@\spl@ttoksb\the\t@@\t@end}%
\def\pr@parec{\e@\spl@ttoksc\the\t@@\t@end}%
\def\pr@pared{\t@@{}\t@@@{}}%

\def\n@metop{\hbox to 0pt{\kern-\d@nkinrightdim\vbox to
\wd0{\vss\offinterlineskip\hsize\ht0\pr@pare
\centerline{\the\t@@@@}\kern2pt\d@nkinball}\hss}}%
\def\n@medown{\hbox to 0pt{\kern-\d@nkinrightdim%
\vtop to \wd0{\offinterlineskip\hsize\ht0\pr@pare%
\d@nkinball\kern2pt\centerline{\the\t@@@@}\vss}\hss}}%
\def\n@meright{\hbox to 0pt{\kern-\d@nkinrightdim\vtop to
\wd0{\offinterlineskip\hsize\ht0\pr@pare
\hbox{\d@nkinball\kern2pt\vbox to 
\wd0{\parindent0pt\vss\rlap{\the\t@@@@}\vss}}\vss}\hss}}%

\def\n@meleft{\hbox to 0pt{\hss\kern-\d@nkinrightdim\vtop to
\wd0{\offinterlineskip\hsize\ht0\pr@pare
\hbox{\vbox to
\wd0{\parindent0pt\vss\llap{\the\t@@@@\kern2pt}\vss}\d@nkinball}\vss}}}

\def\Tl@ne#1{\kern\d@nkinrightdim
\vbox{\hbox to #1{}\hrule}\kern\d@nkinrightdim}%
\def\tl@ne#1{\kern\d@nkinrightdim
\vbox{\hbox to #1{\vrule width1.5mm height.4pt
\if@ne\dimen0 1mm\else\dimen0 2mm\fi
\leaders\hbox to \dimen0{\hss.\hss}\hfill
\vrule width1.5mm height.4pt}}\kern\d@nkinrightdim
\let\hl@ne\Tl@ne}%
\def\b{\global\let\hl@ne\tl@ne}%
\def\vl@ne#1{\raise\dimen0\hbox{\vbox to #1{}\vrule}\kern-.4pt}%
\def\d@nkinone{\dimen9\ht0\advance\dimen9-.4pt
\raise-0.5\dimen9\n@medown}%
\def\d@nkinrightfork{%
\hbox{\kern\d@nkinrightdim\kern-0.1464466094\wd0
\raise-\f@rkraise
\vbox{\offinterlineskip\parindent0pt
\dimen0 0.70710678108\ht0
\if@ne
\advance\dimen0 \ht1
\else
\advance\dimen0 2\ht1
\fi
\hbox{\raise0.85355339059\wd0
\@pline%
\kern0.35355339059\wd0
\raise\dimen0
\n@meright}%
\kern-0.1464466094\ht0
\hbox{\raise0.85355339059\wd0
\d@wnline%
\kern-0.1464466094\wd0
\kern\d@nkinrightdim\n@meright}}}}%
\def\d@nkinleftfork{\kern.5\wd0\raise-\f@rkraise
\vbox{\offinterlineskip\parindent0pt
\dimen0 0.70710678108\ht0
\hbox{\d@wnline
\kern-2\wd1
\if@ne\dimen9 \ht1\else\dimen9 2\ht1\fi
\advance\dimen9-0.1464466094\ht0\if@ne\kern\wd1\fi
\raise\dimen9
\n@meleft}%
\kern0.70710678108\ht0
\hbox{\raise0.85355339059\wd0
\@pline%
\if@ne\kern-\wd1\else\kern-2\wd1\fi
\n@meleft}}\kern0.35355339059\wd0
\if@ne\kern\wd1  \else\kern2\wd1\fi
}%
\def\d@nkinright{\hl@ne\d@nkinstep\d@nkinone}%
\def\d@nkinup{{\vl@ne\d@nkinstep\kern0.5\wd0\raise\dimen2\n@meleft}%
\kern-0.5\wd0}%
\def\d@nkinuptwo{\vl@ne\d@nkinstep\kern0.5\wd0\raise\dimen2\n@meleft%
\kern-0.5\wd0\dimen9\dimen2
\advance\dimen9 0.5\ht0
\raise\dimen9\hbox{\vl@ne\d@nkinstep\kern0.5\wd0\raise\dimen2\n@meleft}%
\kern-0.5\wd0}%
\toksdef\t@0\toksdef\t@@1\toksdef\t@@@2%
\toksdef\t@@@@3%

\def\dynkininit#1{%
\t@{}\t@@{}\t@@@{}\t@@@@{}%
\dimendef\d@nkinstep4
\e@\d@nkinstep\d@nkinsize
\setbox0\hbox{\ifd@dotchoice\dynkinball{#1}\else\dynkinball\fi}%
\dimendef\d@nkinraise6
\dimendef\d@nkinrightdim7
\d@nkinrightdim\wd0
\divide\d@nkinrightdim2
\font\l@net=line10
\dimendef\f@rkraise3
\setbox1\hbox{\l@net\char\d@nkincount}%
\f@rkraise 1.20710678108\ht0
\advance\f@rkraise-.2pt
\if@ne
\advance\f@rkraise \ht1
\else
\advance\f@rkraise 2\ht1 
\fi
\dimen0\ht0
\divide\dimen0 2
\dimen2\d@nkinstep
\advance\dimen2\dimen0
\let\hl@ne\Tl@ne
\edef\@pline{\vbox{\offinterlineskip\parindent0pt\l@net
\if@ne\hsize \wd1\char\d@nkincount
\else
\hsize 2\wd1
\LINE{\hskip\wd1\char\d@nkincount\hss}\par
\char\d@nkincount\fi}}%
\edef\d@wnline{\vbox{\advance\d@nkincount"40\l@net
\offinterlineskip\parindent0pt\if@ne\char\d@nkincount\hsize \wd1
\else
\hsize 2\wd1
\char\d@nkincount\par
\LINE{\hskip\wd1\char\d@nkincount\hss}\fi}}%
}

\def\count@rga#1 #2,{\count255 #1\count253 #1\dynkininit{#2}\next@rg#2,}
\def\count@rgbc#1 #2 {\count255 #1\count253 #1\dynkininit{#2}\next@rg#2 }
\def\count@rgd#1 {\count255 #1\count253 #1\dynkininit{}\next@rg}
\def\init@rga#1,{\dynkininit{#1}\next@rg#1,}
\def\init@rgbc#1 {\dynkininit{#1}\next@rg#1 }
\def\init@rgd{\dynkininit{}\next@rg}

\def\d@loopargsa#1{\def\next@rg##1,##2 {\t@@=\e@{\the\t@@##1\@t##2\@t}\next}
\def\next{\ifnum\count255>1\advance\count255 -1\e@\next@rg\else#1\fi}}
\def\d@loopargsbc#1{\def\next@rg##1 {\t@@=\e@{\the\t@@##1\@t}\next}%
\def\next{\ifnum\count255>1\advance\count255 -1\e@\next@rg\else#1\fi}}
\def\d@loopargsd#1{\def\next@rg{#1}}

\def\read@rgs#1\afterread{
\d@loopargs{#1}
\ifrepe@t
\e@\count@rg
\else
\count255 \count254
\count253 \count254
\e@\init@rg
\fi}

\newif\ifrepe@t
\def\d@nkrepeat{rightrepeat}
\def\d@nkstop{st@p}
\def\d@one{1}
\def\d@right{1}
\def\d@rightrepeat{0}
\def\d@up{1}
\def\d@uptwo{2}
\def\d@rightfork{2}
\def\d@leftfork{2}

\def\countd@nkin#1,{\def\d@rep{#1}%
\ifx\d@nkrepeat\d@rep\repe@ttrue\fi\ifx\d@nkstop\d@rep\else
\advance\count254 \csname d@#1\endcsname
\e@\countd@nkin\fi}

\def\dod@nkin#1,{\def\d@rep{#1}%
\ifx\d@nkstop\d@rep\else
\csname d@nkin#1\endcsname
\e@\dod@nkin\fi}

\def\d@nkinrightrepeat{\loop
\ifnum\count253>0
\advance\count253 -1
\d@nkinright%
\repeat}

\d@setdots
\def\e@ttok#1#2\@nd{\t@@={#2}}

\def\rawdynkin#1;{\count254 0
\repe@tfalse
\countd@nkin#1,st@p,%
\bgroup
\read@rgs
\advance\count253 -\count254
\hbox{
\dod@nkin#1,st@p,%
\kern\d@nkinrightdim}%
\egroup\afterread
}
\def\dynkinA{\rawdynkin one,rightrepeat;}
\def\dynkinD{\rawdynkin one,rightrepeat,rightfork;}
\def\dynkinDtilde#1 {\count252 #1\advance\count 252 1
\rawdynkin leftfork,one,rightrepeat,rightfork;{\count252} }
\def\dynkinE{\rawdynkin one,right,right,up,rightrepeat;}
\def\dynkinEtilde#1 {\let\d@it\relax\ifcase#1\or\or\or\or\or\or
\def\d@it{\rawdynkin one,right,right,uptwo,right,right;}\or
\def\d@it{\rawdynkin one,right,right,right,up,right,right,right;}\else
\count252 #1\advance\count 252 1
\def\d@it{\rawdynkin one,right,right,up,rightrepeat;{\count252} }\fi
\d@it}

\def\d@nkinbullet{\vbox to 3.85pt{\kern-.1pt\hbox to 3.7pt
{\kern-.43pt$\bullet$\hss}\vss}}

\def\d@nkinbulletandcirc#1{\vbox to 3.85pt{\kern-.1pt\hbox to 3.7pt
{\kern-.43pt\ifx #1c$\circ$\else$\bullet$\fi\hss}\vss}}

\def\d@nkinbigcirc{\vbox to 8.4pt{\kern-.1pt\hbox to 8.8pt
{\kern-.4pt$\bigcirc$\hss}\vss}}

\setdynkin ball=bullet
\catcode`\@=12
